\documentclass[aap,seceqn,nameyear,MSNbibl,dvips]{arximspdf}
\usepackage{mathrsfs}

\doi{10.1214/09-AAP646}
\volume{20}
\issue{3}
\pubyear{2010}
\firstpage{1005}
\lastpage{1028}

\makeatletter

\newtheorem{theorem}{Theorem}[section]
\newtheorem{proposition}{Proposition}[section]
\newtheorem{lemma}{Lemma}[section]

\newproclaim{definition}{Definition}[section]

\makeatother

\begin{document}
\begin{frontmatter}

\title{An asymptotic sampling formula for the coalescent
with Recombination}
\runtitle{Asymptotic sampling formula for a two-locus model}

\begin{aug}
\author[A]{\fnms{Paul A.} \snm{Jenkins}\thanksref{t1}\ead[label=e1]{pauljenk@eecs.berkeley.edu}} and
\author[B]{\fnms{Yun S.} \snm{Song}\corref{}\thanksref{t1,t2}\ead[label=e2]{yss@stat.berkeley.edu}}
\runauthor{P. A. Jenkins and Y. S. Song}
\affiliation{University of California, Berkeley}
\address[A]{Computer Science Division\\
University of California, Berkeley\\
Berkeley, California 94720\\
USA\\
\printead{e1}}
\address[B]{Department of Statistics\\
\quad and Computer Science Division\\
University of California, Berkeley\\
Berkeley, California 94720\\
USA\\
\printead{e2}}
\end{aug}

\thankstext{t1}{Supported in part by an NIH Grant R00-GM080099.}
\thankstext{t2}{Supported in part by an Alfred P. Sloan Research
Fellowship and a Packard Fellowship for
Science and Engineering.}

\received{\smonth{5} \syear{2009}}
\revised{\smonth{9} \syear{2009}}

%
\begin{abstract}
Ewens sampling formula (ESF) is a one-parameter family of probability
distributions with a number of intriguing combinatorial connections.
This elegant closed-form formula first arose in biology as the
stationary probability distribution of a sample configuration at one
locus under the infinite-alleles model of mutation. Since its discovery
in the early 1970s, the ESF has been used in various biological
applications, and has sparked several interesting mathematical
generalizations. In the population genetics community, extending the
underlying random-mating model to include recombination has received
much attention in the past, but no general closed-form sampling formula
is currently known even for the simplest extension, that is, a model
with two loci. In this paper, we show that it is possible to obtain
useful closed-form results in the case the population-scaled
recombination rate $\rho$ is large but not necessarily infinite.
Specifically, we consider an asymptotic expansion of the two-locus
sampling formula in inverse powers of $\rho$ and obtain closed-form
expressions for the first few terms in the expansion. Our asymptotic
sampling formula applies to arbitrary sample sizes and configurations.
\end{abstract}

%

\end{frontmatter}

\section{Introduction}
\label{int}
The probability of a sample configuration provides a useful ground for
analyzing genetic data. Popular applications include obtaining maximum
likelihood estimates of model parameters and performing ancestral
inference [see \citet{ste2001}]. In principle, model-based
full-likelihood analyses, such as that based on the coalescent
[Kingman (\citeyear{kin1982b}, \citeyear{kin1982a})], should be among
the most powerful
methods since they make full use of the data. However, in most cases,
it is intractable to obtain a closed-form formula for the probability
of a given data set. A~well-known exception to this hurdle is the Ewens
sampling formula (ESF), which describes the stationary probability
distribution of a sample configuration under the one-locus
infinite-alleles model in the diffusion limit [\citet{ewe1972}].
Notable biological applications of this closed-form formula include the
test of selective neutrality [see
\citet{Watterson1977}, Slatkin (\citeyear{SlatkinGenetRes94}, \citeyear{SlatkinGenetRes96})].
\citet
{hop1984}
provided a P\'{o}lya-like urn model interpretation of the formula, and
recently \citet{Griffiths2005p1491} provided a new combinatorial
proof of the ESF and extended the framework to obtain new results for
the case with a variable population size. We refer the reader to the
latter paper for a nice summary of previous works related to the ESF.
Note that the ESF also arises in several interesting contexts outside
biology, including random partition structures and Bayesian statistics;
see \citet{ABTbook} for examples of intricate combinatorial
connections. The ESF is a special case of the two-parameter sampling
formula constructed by Pitman (\citeyear{pitman1992}, \citeyear
{pitman1995}) for
exchangeable random partitions.

\citet{gol1984} considered generalizing the infinite-alleles
model to include recombination and constructed a recursion relation
satisfied by the two-locus sampling probability distribution at
stationarity in the diffusion limit. \citet{ethgri1990} later
undertook a more mathematical analysis of the two-locus model and
provided several interesting results. However, to date, a general
closed-form formula for the two-locus sampling distribution remains
unknown. Indeed, it is widely recognized that recombination adds a
formidably challenging layer of complexity to population genetics
analysis. Because obtaining exact analytic results in the presence of
recombination is difficult, recent research has focused on developing
sophisticated and computationally-intensive Monte Carlo techniques.
Examples of such techniques applied to the coalescent include Monte
Carlo simulations [see Hudson (\citeyear{hud1985}, \citeyear{hud2001})],
importance sampling [see \citet{grimar1996},
\citet{stedon2000}, \citet{feadon2001},
De Iorio and Griffiths (\citeyear{deigri2004a}, \citeyear{deigri2004b}),
\citet{grietal2008}] and Markov
chain Monte Carlo methods [see \citet{kuhetal2000},
\citet{nie2000}, \citet{wanran2008}].

Being the simplest model with recombination, the two-locus case has
been extensively studied in the past [\citet{gri1981},
\citet{gol1984}, \citet{hud1985}, \citet{ethgri1990},
\citet{gri1991}] and a renewed wave of interest was
recently sparked by \citet{hud2001}, who proposed a composite
likelihood method which uses two-locus sampling probabilities as
building blocks. LDhat, a widely-used software package for estimating
recombination rates, is based on this composite likelihood approach,
and it has been used to produce a fine-scale map of recombination rate
variation in the human genome
[\citet{McVeanScience04}, \citet{MyersScience05}]. LDhat
assumes a
symmetric diallelic recurrent mutation model at each locus and relies
on the importance sampling scheme proposed by \citet{feadon2001}
for the coalescent with recombination, to generate exhaustive lookup
tables containing two-locus probabilities for all inequivalent sample
configurations and a range of relevant parameter values. This process
of generating exhaustive lookup tables is very computationally
expensive. A fast and accurate method of estimating two-locus
probabilities would be of practical value.

In this paper, we revisit the tantalizing open question of whether a
closed-form sampling formula can be found for the coalescent with
recombination. We show that, at least for the two-locus
infinite-alleles model with the population-scaled recombination rate
$\rho$ large but not necessarily infinite, it is possible to obtain
useful closed-form analytic results. Note that the aforementioned Monte
Carlo methods generally become less efficient as $\rho$ increases.
Those methods involve sampling a large collection of genealogical
histories consistent with the observed sample configuration, and, when
$\rho$ is large, the sampled genealogies tend to be very complicated;
they typically contain many recombination events, and it may take a
long time for every locus to reach a most recent common ancestor.
However, contrary to this increased complexity in the standard
coalescent, we actually expect the evolutionary dynamics to be easier
to describe for large $\rho$, since the loci under consideration would
then be less dependent. Hence, it seems reasonable to conjecture that
there may exist a stochastic process simpler than the standard
coalescent with recombination that describes the relevant degrees of
freedom in the large $\rho$ limit. We believe that our sampling formula
may provide some hints as to what that dual process should be.

The work discussed here generalizes previous results [\citet
{gol1984}, \citet{ethgri1990}] for $\rho=\infty$, in which case
the loci
become independent and the two-locus sampling distribution is given by
a product of one-locus ESFs. Our main results can be summarized as
follows.

\subsection*{Main results} Consider the diffusion limit of the two-locus
infinite-alleles model with population-scaled mutation rates $\theta_A$
and $\theta_B$ at the two loci. For a sample configuration
${\mathbf{n}}$ (defined later in the text), we use
$q({\mathbf{n}}\mid\theta_A,\theta_B,\rho)$ to denote the probability of
observing ${\mathbf{n}}$ given the parameters $\theta_A,\theta_B$ and
$\rho$. For an arbitrary ${\mathbf{n}}$, our goal is to find an
asymptotic expansion of $q({\mathbf{n}}\mid\theta_A,\theta_B,\rho)$ in
inverse powers of~$\rho$, that is, for large values of the recombination
rate $\rho$, our goal is to find
%
\[
q({\mathbf{n}}\mid\theta_A,\theta_B,\rho) = q_0({\mathbf{n}}\mid\theta
_A,\theta_B) + \frac{q_1({\mathbf{n}}\mid\theta_A
,\theta_B)}{\rho} +
\frac{q_2({\mathbf{n}}\mid\theta_A,\theta_B)}{\rho^2} + O
\biggl(\frac{1}{\rho^3} \biggr),
\]
%
where $q_0, q_1$, and $q_2$ are independent of $\rho$. As mentioned
before, $q_0({\mathbf{n}}\mid\theta_A,\theta_B)$ is given by a product of
one-locus ESFs. In this paper, we derive a closed-form formula for the
first-order term $q_1({\mathbf{n}}\mid\theta_A,\theta_B)$. Further, we
show that the second-order term $q_2({\mathbf{n}}\mid\theta_A,\theta_B)$
can be decomposed into two parts, one for which we obtain a closed-form
formula and the other that satisfies a simple strict recursion. The
latter can be easily evaluated using dynamic programming. Details of
these results are described in Section \ref{sec:asy}. In a similar
vein, in Section \ref{sec:joi}, we obtain a simple asymptotic formula
for the joint probability distribution of the number of alleles
observed at the two loci.

We remark that our work has practical value in genetic analysis. While
this paper was under review, we applied the technique developed here to
obtain analogous results for an arbitrary finite-alleles recurrent
mutation model. See \citet{jenson2009} for details. In that
paper, we performed an extensive assessment of the accuracy of our
results for a particular finite-alleles model of mutation, and showed
that they may be accurate even for moderate values of $\rho$, including
a range that is of biological interest. The accuracy (not discussed
here) of our results for the infinite-alleles model is very similar to
that finite-alleles case.

\section{Preliminaries}
\label{sec:pre}

In this section, we review the ESF for the one-locus infinite-alleles
model, as well as Golding's (\citeyear{gol1984}) recursion relation
for the two-locus generalization. Our notational convention generally
follows that of \citet{ethgri1990}.

Given a positive integer $k$, $[k]$ denotes the $k$-set
$\{1,\ldots,k\}$. For a nonnegative real number $x$ and a positive
integer $n$, $(x)_{n} := x(x+1)\cdots(x+n-1)$ denotes the $n$th
ascending factorial of $x$. We use $\mathbf{0}$ to denote either a
vector or a matrix of all zeroes; it will be clear from context which
is intended. Throughout, we consider the diffusion limit of a neutral
haploid exchangeable model of random mating with constant population
size $2N$. We refer to the haploid individuals in the population as
gametes.

\subsection{Ewens sampling formula for the one-locus model}
\label{sec:ESF}
In the one-locus model, a sample configuration is denoted by a vector
of multiplicities ${\mathbf{n}}= (n_1,\ldots,n_K)$, where $n_i$ denotes
the number of gametes with allele $i$ at the locus and $K$ denotes the
total number of distinct allelic types observed. We use $n$ to denote
$\sum_{i=1}^K n_i$, the total sample size. Under the infinite-alleles
model, any two gametes can be compared to determine whether or not they
have the same allele, but it is not possible to determine how the
alleles are related when they are different. Therefore, allelic label
is arbitrary. The probability of a mutation event at the locus per
gamete per generation is denoted by $u$. In the diffusion limit,
$N\to\infty$ and $u\to0$ with the population-scaled mutation rate
$\theta= 4N u$ held fixed. Each mutation gives rise to a new allele
that has never been seen before in the population. For the one-locus
model just described, \citet{ewe1972} obtained the following
result.
\begin{proposition}[(Ewens)]
At stationarity in the diffusion limit of the one-locus
infinite-alleles model with the scaled mutation parameter $\theta$, the
probability of an unordered sample configuration
${\mathbf{n}}=(n_1,\ldots,n_K)$ is given by
%
%
\begin{equation}\label{ESF}
p({\mathbf{n}}\mid\theta) = \frac{n!}{n_1\cdots
n_K}\frac{1}{\alpha_1!\cdots\alpha_n!}\frac{\theta^K}{(\theta)_{n}},
\end{equation}
%
where $\alpha_i$ denotes the number of allele types represented $i$
times, that is, $\alpha_i:=| \{ k \mid n_k = i\}|$.
\end{proposition}

Let $\mathscr{A}_n$ denote an \textit{ordered} configuration of $n$
sequentially sampled gametes such that the corresponding unordered
configuration is given by ${\mathbf{n}}$. By exchangeability, the
probability of $\mathscr{A}_n$ is invariant under all permutations of
the sampling order. Hence, we can write this probability of an ordered
sample as $q({\mathbf{n}})$ without ambiguity. It is given by
%
%
\begin{equation}
\label{q(A)}\quad
q({\mathbf{n}}\mid\theta)
= p({\mathbf{n}}\mid\theta) \biggl[\frac{n!}{\prod_{i=1}^K
n_i!}\frac
{1}{\alpha_1!\cdots\alpha_n!} \biggr]^{-1}
= \Biggl[ \prod_{i=1}^K (n_i - 1)! \Biggr]\frac{\theta^K}{(\theta)_{n}},
\end{equation}
%
which follows from the fact that there are $\frac{n!}{\prod_{i=1}^K
n_i!} \frac{1}{\alpha_1!\cdots\alpha_n!}$ orderings corresponding to
${\mathbf{n}}$ [\citet{hop1984}]. To understand what we mean by
ordered and unordered sampling configurations, it is helpful to relate
the Ewens sampling formula to the theory of random partitions. If the
gametes are labeled in order of appearance by $1,\ldots,n$, then the
resulting sample configuration defines a random partition of $[n]$,
with gametes belonging to the same block if and only if they have the
same allele. The quantity $q({\mathbf{n}}\mid\theta)$ is then the
probability of a \textit{particular} partition of $[n]$ whose block sizes
are given by the entries in ${\mathbf{n}}$, while the quantity
$p({\mathbf{n}}\mid\theta)$ is the probability of observing \textit{any}
partition of $[n]$ with these block sizes. For example, if
${\mathbf{n}}= (2,1,1)$, then there are six partitions of $[4]$ with
these block sizes, and so $p({\mathbf{n}}\mid\theta) =
6q({\mathbf{n}}\mid\theta)$. It is often more convenient to work with
an ordered sample than with an unordered sample. In this paper, we will
work with the former; that is, we will work with
$q({\mathbf{n}}\mid\theta)$ rather than $p({\mathbf{n}}\mid\theta)$.

In the coalescent process going backward in time, at each event a
lineage is lost either by coalescence or mutation. By consideration of
the most recent event back in time, one can show that
$q({\mathbf{n}}\mid\theta)$ satisfies
%
%
\begin{eqnarray}\label{one-locus-recursion}
n(n-1+\theta)q({\mathbf{n}}\mid\theta) &=& \sum_{i=1}^K n_i(n_i -
1)q({\mathbf{n}}-{\mathbf{e}}_i\mid\theta)\nonumber\\[-8pt]\\[-8pt]
&&{}+
\theta\sum_{i=1}^K \delta_{n_i,1} q({\mathbf{n}}-{\mathbf
{e}}_i\mid\theta),\nonumber
\end{eqnarray}
%
where $\delta_{n_i,1}$ is the Kronecker delta and ${\mathbf{e}}_i$ is a
unit vector with the $i$th entry equal to one and all other entries
equal to zero. The boundary condition is $q({\mathbf{e}}_i\mid\theta) =
1$ for all $i\in[K]$, and $q({\mathbf{n}}\mid\theta)$ is defined to be
zero if ${\mathbf{n}}$ contains any negative component. It can be
easily verified that the formula of $q({\mathbf{n}}\mid\theta)$ shown
in (\ref{q(A)}) satisfies the recursion (\ref{one-locus-recursion}).

\citet{ewe1972} also obtained the following result regarding the
number of allelic types.
%
\begin{proposition}[(Ewens)]
Let $K_n$ denote the number of distinct allelic types observed in a
sample of size $n$. Then
%
%
\begin{equation}\label{eq:K_n}
\mathbb{P}(K_n = k\mid\theta) = \frac{s(n,k)\theta^k}{(\theta)_{n}},
\end{equation}
where $s(n,k)$ are the unsigned Stirling numbers of the first kind.
Note that $(\theta)_{n} = s(n,1)\theta+ s(n,2)\theta^2 + \cdots+
s(n,n)\theta^n$.
\end{proposition}

It follows from (\ref{ESF}) and (\ref{eq:K_n}) that
$K_n$ is a sufficient statistic for $\theta$.

\subsection{Golding's recursion for the two-locus case}

Golding (\citeyear{gol1984}) first generalized the one-locus recursion
(\ref{one-locus-recursion}) to two loci, and \citet{ethgri1990}
later undertook a more mathematical study of the model. We denote the
two loci by $A$ and $B$, and use $\theta_A$ and $\theta_B$ to denote
the respective population-scaled mutation rates. We use $K$ and $L$ to
denote the number of distinct allelic types observed at locus $A$ and
locus $B$, respectively. The population-scaled recombination rate is
denoted by $\rho= 4Nr$, where $r$ is the probability of a recombination
event between the two loci per gamete per generation. A key observation
is that to obtain a closed system of equations, the type space must be
extended to allow some gametes to be specified only at one of the two
loci.
\begin{definition}[(Extended sample configuration for two loci)]
The two-locus sample configuration is denoted by ${\mathbf{n}}=
({\mathbf{a}},{\mathbf{b}},{\mathbf{c}})$, where
${\mathbf{a}}=(a_1,\ldots,a_K)$ with $a_i$ being the number of gametes
with allele $i$ at locus $A$ and unspecified alleles at locus $B$,
${\mathbf{b}}=(b_1,\ldots,b_L)$ with $b_j$ being the number of gametes
with unspecified alleles at locus $A$ and allele $j$ at locus $B$, and
${\mathbf{c}}=(c_{ij})$ is a $K \times L$ matrix with $c_{ij}$ being
the multiplicity of gametes with allele $i$ at locus $A$ and allele $j$
at locus $B$. Further, we define
\begin{eqnarray*}
a &=& \sum_{i=1}^K a_i,\qquad
c_{i\cdot} = \sum_{j=1}^L c_{ij},\qquad
c = \sum_{i=1}^K \sum_{j=1}^L c_{ij}, \\
b &=& \sum_{j=1}^L b_j,\qquad
c_{\cdot j} = \sum_{i=1}^K c_{ij},\qquad
n = a + b + c.
\end{eqnarray*}
\end{definition}

We use $q({\mathbf{a}},{\mathbf{b}},{\mathbf{c}})$ to denote the
sampling probability of an ordered sample with configuration
$({\mathbf{a}},{\mathbf{b}},{\mathbf{c}})$. For ease of notation, we do
not show the dependence on parameters. For $0 \leq\rho< \infty$,
Golding's (\citeyear{gol1984}) recursion for
$q({\mathbf{a}},{\mathbf{b}},{\mathbf{c}})$ takes the following form:
%
%
\begin{eqnarray}\label{golding}
&&[n(n-1) + \theta_A(a+c) + \theta_B(b+c) + \rho
c]q({\mathbf
{a}},{\mathbf{b}},{\mathbf{c}}) \nonumber\\
& &\qquad= \sum_{i=1}^K a_i(a_i-1+2c_{i\cdot})q({\mathbf{a}}-{\mathbf{e}}
_i,{\mathbf{b}},{\mathbf{c}})\nonumber\\
&&\qquad\quad{} + \sum_{j=1}^L b_j(b_j-1+2c_{\cdot
j})q({\mathbf{a}},{\mathbf{b}}-{\mathbf{e}}_j,{\mathbf{c}})
\nonumber\\
& &\qquad\quad{} + \sum_{i=1}^K\sum_{j=1}^L [
c_{ij}(c_{ij}-1)q({\mathbf{a}},{\mathbf{b}},{\mathbf{c}}-{\mathbf{e}}_{ij})
\nonumber\\
&&\qquad\quad\hspace*{45.2pt}{} + 2 a_i b_j q({\mathbf{a}}-{\mathbf{e}}_i,{\mathbf{b}}-{\mathbf
{e}}_j,{\mathbf
{c}}+{\mathbf{e}}_{ij}) ] \\
& &\qquad\quad{} + \theta_A\sum_{i=1}^K \Biggl[\sum_{j=1}^L \delta
_{a_i+c_{i\cdot
},1}\delta_{c_{ij},1} q({\mathbf{a}},{\mathbf{b}}+{\mathbf
{e}}_j,{\mathbf
{c}}-{\mathbf{e}}_{ij})\nonumber\\
&&\qquad\quad\hspace*{92.8pt}{} + \delta_{a_i,1}\delta_{c_{i\cdot},0}
q({\mathbf{a}}-
{\mathbf{e}}_i,{\mathbf{b}},{\mathbf{c}}) \Biggr] \nonumber\\
& &\qquad\quad{} + \theta_B\sum_{j=1}^L \Biggl[\sum_{i=1}^K \delta_{b_j +
c_{\cdot j},1}
\delta_{c_{ij},1} q({\mathbf{a}}+{\mathbf{e}}_i,{\mathbf
{b}},{\mathbf
{c}}-{\mathbf{e}}_{ij})\nonumber\\
&&\qquad\quad\hspace*{90.1pt}{} +\delta_{b_j,1}\delta_{c_{\cdot j},0}
q({\mathbf
{a}},{\mathbf{b}}- {\mathbf{e}}_j,{\mathbf{c}}) \Biggr] \nonumber\\
& &\qquad\quad{} + \rho\sum_{i=1}^K\sum_{j=1}^L c_{ij}q({\mathbf{a}}+{\mathbf{e}}
_i,{\mathbf{b}}+{\mathbf{e}}_j,{\mathbf{c}}-{\mathbf
{e}}_{ij}).\nonumber
\end{eqnarray}
%
Relevant boundary conditions are $q({\mathbf{e}}_i,\mathbf{0},\mathbf
{0}) = q(\mathbf{0},{\mathbf{e}}_j,\mathbf{0}) = 1$ for all $i\in[K]$
and $j\in[L]$. For notational convenience, we deviate from
\citet{ethgri1990} and allow each summation to range over all
allelic types. To be consistent, we define
$q({\mathbf{a}},{\mathbf{b}},{\mathbf{c}}) = 0$ whenever any entry in
${\mathbf{a}}$, ${\mathbf{b}}$ or ${\mathbf{c}}$ is negative.

For ease of discussion, we define the following terms.
\begin{definition}[(Degree)]
The degree of $q({\mathbf{a}},{\mathbf{b}},{\mathbf{c}})$ is defined
to be $a+b+2c$.
\end{definition}
\begin{definition}[(Strictly recursive)]
We say that a recursion relation is strictly recursive if it contains
only a single term of the highest degree.
\end{definition}

Except in the special case $\rho= \infty$, a closed-form solution for
$q({\mathbf{a}},{\mathbf{b}},{\mathbf{c}})$ is not known. Notice that
the terms $q({\mathbf{a}}-{\mathbf{e}}_i,{\mathbf{b}}-{\mathbf{e}}
_j,{\mathbf{c}}+{\mathbf{e}}_{ij})$ and
$q({\mathbf{a}}+{\mathbf{e}}_i,{\mathbf{b}}+{\mathbf{e}}_j,{\mathbf
{c}}-{\mathbf{e}}_{ij})$ on the right-hand side of (\ref{golding}) have
the same degree as $q({\mathbf{a}},{\mathbf{b}},{\mathbf{c}})$ on the
left-hand side. Therefore, (\ref{golding}) is not strictly recursive.
For each degree, we therefore need to solve a system of coupled
equations, and this system grows very rapidly with~$n$. For example,
for a sample with $a=0, b=0$ and $c=40$, computing
$q(\mathbf{0},\mathbf{0},{\mathbf{c}})$ requires solving a system of
more than 20,000 coupled equations [\citet{hud2001}]; this is
around the limit of sample sizes that can be handled in a reasonable
time. In the following section, we revisit the problem of obtaining a
closed-form formula for $q({\mathbf{a}},{\mathbf{b}},{\mathbf{c}})$ and
obtain an asymptotic expansion for large $\rho$.

\section{An asymptotic sampling formula for the two-locus case}
\label{sec:asy}

For large $\rho$, our objective is to find an asymptotic expansion of
the form
%
%
\begin{equation}
\label{expansion}
q({\mathbf{a}},{\mathbf{b}},{\mathbf{c}}) = q_0({\mathbf
{a}},{\mathbf{b}},{\mathbf{c}}) +
\frac{q_1({\mathbf{a}},{\mathbf{b}},{\mathbf{c}})}{\rho} + \frac
{q_2({\mathbf{a}},{\mathbf{b}},{\mathbf{c}})}{\rho^2} +
O \biggl(\frac{1}{\rho^{3}} \biggr),
\end{equation}
where $q_0,q_1$ and $q_2$ are independent of $\rho$.
Our closed-form formulas will be expressed using
the following notation.
\begin{definition}
For a given multiplicity vector ${\mathbf{a}}=(a_1,\ldots,a_K)$ with
$a=\sum_{i=1}^K a_i$, we define
%
%
\begin{equation}
\label{ESFa}
q^A({\mathbf{a}}) = \Biggl[\prod_{i=1}^K (a_i - 1)! \Biggr] \frac
{\theta_A^K}{(\theta_A)_{a}}.
\end{equation}
Similarly, for a given multiplicity vector ${\mathbf{b}}=(b_1,\ldots
,b_L)$ with $b=\sum_{i=1}^L b_i$,
we define
%
%
\begin{equation}
\label{ESFb}
q^B({\mathbf{b}}) = \Biggl[\prod_{j=1}^L (b_j - 1)! \Biggr]\frac
{\theta_B^L}{(\theta_B)_{b}}.
\end{equation}
As discussed in Section \ref{sec:ESF}, $q^A$ (respectively, $q^B$)
gives the probability of an ordered sample taken from locus $A$
(respectively, $B$).
\end{definition}
\begin{definition}[(Marginal configuration)]
We use ${\mathbf{c}}_A= (c_{i\cdot})_{i\in[K]}$ and ${\mathbf{c}}_B=
(c_{\cdot j})_{j\in[L]}$ to denote the marginal sample configurations
of ${\mathbf{c}}$ restricted to locus $A$ and locus $B$, respectively.
\end{definition}

The leading-order term $q_0({\mathbf{a}},{\mathbf{b}},{\mathbf{c}})$ is
equal to $q({\mathbf{a}},{\mathbf{b}},{\mathbf{c}})$ when $\rho=
\infty$, in which case the two loci are independent. Theorem 2.3 of
\citet{ethgri1990} states that
$q_0(\mathbf{0},\mathbf{0},{\mathbf{c}})=q^A({\mathbf
{c}}_A)q^B({\mathbf{c}}_B)$. More generally, one can obtain the
following result for the leading-order contribution.
%
%
\begin{proposition}
\label{prop:leading}
In the asymptotic expansion (\ref{expansion}) of the two-locus
sampling formula, the zeroth-order term $q_0({\mathbf{a}},{\mathbf
{b}},{\mathbf{c}})$
is given by
%
%
\begin{equation}
\label{leading}
q_0({\mathbf{a}},{\mathbf{b}},{\mathbf{c}}) = q^A({\mathbf
{a}}+{\mathbf{c}}_A)q^B({\mathbf{b}}+{\mathbf{c}}_B).
\end{equation}
\end{proposition}
%

Although this result is intuitively obvious, in Section \ref
{sec:proof_leading} we provide
a detailed new proof, since it well illustrates our general strategy.
One of the main results of this paper is
a closed-form formula for the next order term $q_1({\mathbf
{a}},{\mathbf{b}},{\mathbf{c}})$.
The case with ${\mathbf{c}}= \mathbf{0}$ admits a particularly simple
solution.
\begin{lemma}
\label{lemma:first}
In the asymptotic expansion (\ref{expansion}) of the two-locus
sampling formula, the first-order term satisfies
\[
q_1({\mathbf{a}},{\mathbf{b}},\mathbf{0}) = 0
\]
for arbitrary ${\mathbf{a}}$ and ${\mathbf{b}}$.
\end{lemma}

That $q_1({\mathbf{a}},{\mathbf{b}},\mathbf{0})$ vanishes is not
expected a priori. Below we shall see that $q_2({\mathbf{a}},{\mathbf
{b}},\mathbf{0})
\neq0$ in general.
For an arbitrary configuration matrix ${\mathbf{c}}$ of nonnegative integers,
we obtain the following closed-form formula for $q_1({\mathbf
{a}},{\mathbf{b}},{\mathbf{c}})$.
%
\begin{theorem}
\label{thm:first}
In the asymptotic expansion (\ref{expansion}) of the two-locus
sampling formula, the first-order term
$q_1({\mathbf{a}},{\mathbf{b}},{\mathbf{c}})$ is given by
%
%
\begin{eqnarray}
\label{firstorder}
q_1({\mathbf{a}},{\mathbf{b}},{\mathbf{c}}) &=& \pmatrix{c\cr2} q^A
({\mathbf{a}}+{\mathbf{c}}_A)q^B({\mathbf{b}}+{\mathbf{c}}_B)
\nonumber\\
& &{}
-q^B({\mathbf{b}}+{\mathbf{c}}_B)\sum_{i=1}^K \pmatrix{c_{i\cdot}
\cr2} q^A({\mathbf{a}}+{\mathbf{c}}_A-{\mathbf{e}}_i) \nonumber
\nonumber\\[-8pt]\\[-8pt]
& &{}
-q^A({\mathbf{a}}+{\mathbf{c}}_A)\sum_{j=1}^L \pmatrix{c_{\cdot j}
\cr2} q^B({\mathbf{b}}+{\mathbf{c}}_B-{\mathbf{e}}_j) \nonumber
\nonumber\\
& &{}
+ \sum_{i=1}^K \sum_{j=1}^L \pmatrix{c_{ij} \cr2} q^A({\mathbf
{a}}+{\mathbf{c}}_A-{\mathbf{e}}_i)q^B({\mathbf{b}}+{\mathbf{c}}_B
-{\mathbf{e}}_j)\nonumber
\end{eqnarray}
for arbitrary configurations ${\mathbf{a}},{\mathbf{b}},{\mathbf
{c}}$ of nonnegative integers.
\end{theorem}
%

Lemma \ref{lemma:first} is used in proving Theorem \ref{thm:first}.
A proof of Theorem \ref{thm:first} is provided in Section \ref
{sec:proof_first},
while a proof of Lemma \ref{lemma:first} is given in
Section \ref{sec:proof_lemma_first}.
Note that the functional form of $q_0({\mathbf{a}},{\mathbf
{b}},{\mathbf{c}})$ and
$q_1({\mathbf{a}},{\mathbf{b}},{\mathbf{c}})$ in (\ref{leading}) and
(\ref{firstorder})
has no explicit dependence on mutation;
that is, the dependence on mutation is completely absorbed
into the marginal one-locus probabilities.
It turns out that (\ref{leading}) and (\ref{firstorder})
are \textit{universal} in that they
also apply to an arbitrary finite-alleles model of mutation, with
$q^A$ and $q^B$ replaced with
appropriate marginal one-locus probabilities for the assumed
mutation model. See \citet{jenson2009} for details.

In principle, similar arguments can be used to find the
$(j+1)$th-order term given the $j$th, although a general expression
does not seem to be easy to obtain.
In Section \ref{sec:proof_second}, we provide a proof of the following result
for $q_2({\mathbf{a}},{\mathbf{b}},{\mathbf{c}})$.
%
\begin{theorem}
\label{thm:second}
In the asymptotic expansion (\ref{expansion}) of the two-locus
sampling formula,
the second-order term $q_2({\mathbf{a}},{\mathbf{b}},{\mathbf{c}})$
is of the form
%
%
\begin{equation}
\label{q2concise}
q_2({\mathbf{a}},{\mathbf{b}},{\mathbf{c}}) = q_2({\mathbf
{a}}+{\mathbf{c}}_A,{\mathbf{b}}+{\mathbf{c}}_B,\mathbf{0}) +
\sigma({\mathbf{a}},{\mathbf{b}},{\mathbf{c}}),
\end{equation}
where $\sigma({\mathbf{a}},{\mathbf{b}},{\mathbf{c}})$ is
given by the analytic formula shown in the \hyperref[appendixA]{Appendix},
and $q_2({\mathbf{a}},{\mathbf{b}},\mathbf{0})$ satisfies the following
strict recursion:
%
%
\begin{eqnarray}
\label{q2(a,b,0)}\quad
&&[a(a+\theta_A-1) + b(b+\theta_B-1)]q_2({\mathbf
{a}},{\mathbf
{b}},\mathbf{0}) \nonumber\\
& &\qquad= \sum_{i=1}^K a_i(a_i-1)q_2({\mathbf{a}}-{\mathbf
{e}}_i,{\mathbf
{b}},\mathbf{0}) + \sum_{j=1}^L b_j(b_j-1)q_2({\mathbf{a}},{\mathbf
{b}}-{\mathbf{e}}_j,\mathbf{0}) \nonumber\\
& &\qquad\quad{} + \theta_A\sum_{i=1}^K \delta_{a_i,1}q_2({\mathbf
{a}}-{\mathbf{e}}
_i,{\mathbf{b}},\mathbf{0}) + \theta_B\sum_{j=1}^L\delta
_{b_j,1}q_2({\mathbf{a}},{\mathbf{b}}-{\mathbf{e}}_j,\mathbf{0})
\\
& &\qquad\quad{} + 4
\Biggl[a\theta_A - (\theta_A+ a - 1)\sum
_{i=1}^K \delta
_{a_i,1} \Biggr]
\nonumber\\
&&\qquad\quad\hspace*{15.7pt}{}\times
\Biggl[b\theta_B - (\theta_B+ b - 1)\sum
_{j=1}^L \delta
_{b_j,1} \Biggr]
q^A({\mathbf{a}})q^B({\mathbf{b}})\nonumber
\end{eqnarray}
with boundary conditions $q_2({\mathbf{e}}_i,\mathbf{0},\mathbf{0}) =
q_2(\mathbf{0},{\mathbf{e}}_j,\mathbf{0}) = 0$ for all $i\in[K]$
and $j\in[L]$.
\end{theorem}
%

In contrast to $q_1({\mathbf{a}},{\mathbf{b}},\mathbf{0})$ (cf.
Lemma \ref{lemma:first}),
it turns out that $q_2({\mathbf{a}},{\mathbf{b}},\mathbf{0})$ does not
vanish in general.
We do not have an analytic solution for
$q_2({\mathbf{a}},{\mathbf{b}},\mathbf{0})$, but note that
(\ref{q2(a,b,0)}) is strictly recursive and that
it can be easily solved numerically using
dynamic programming.
Numerical study (not shown) suggests that the relative contribution of
$q_2({\mathbf{a}}+{\mathbf{c}}_A,{\mathbf{b}}+{\mathbf{c}}_B
,\mathbf{0})$ to $q({\mathbf{a}},{\mathbf{b}},{\mathbf{c}})$ is
in most cases extremely small. [See \citet{jenson2009} for
details.]
Deriving an analytic expression for
$\sigma({\mathbf{a}},{\mathbf{b}},{\mathbf{c}})$ in (\ref{q2concise})
is a laborious
task, as the long equation in the \hyperref[appendixA]{Appendix}
suggests. We have
written a computer program to
verify numerically that our analytic result is correct.

\section{Joint distribution of the number of alleles at the two loci
in a sample}
\label{sec:joi}
Following the same strategy as in the previous section, we can obtain the
asymptotic behavior of the joint distribution of the number of
alleles observed at the two loci
in a sample. To make explicit the dependence of these numbers
on the sample size, write the number of alleles at locus $A$ as
$K_{a,b,c}$ and the number of alleles at locus $B$ as
$L_{a,b,c}$. \citet{ethgri1990} proved
that the probability $p(a,b,c;k,l) := \mathbb{P}(K_{a,b,c} = k,
L_{a,b,c} =l)$
satisfies the recursion
%
%
\begin{eqnarray}
\label{KL}\hspace*{25pt}
&&[n(n-1) + \theta_A(a+c) + \theta_B(b+c) + \rho
c]p(a,b,c;k,l)
\nonumber\\
& &\qquad = a(a-1+2c)p(a-1,b,c;k,l) + b(b-1+2c)p(a,b-1,c;k,l) \nonumber\\
& &\qquad\quad {}+ c(c-1)p(a,b,c-1;k,l) + 2a b p(a-1,b-1,c+1;k,l) \nonumber\\
& &\qquad\quad {}+ \theta_A [a p(a-1,b,c;k-1,l) + c p(a,b+1,c-1;k-1,l) ]
\\
& &\qquad\quad {}+ \theta_B [b p(a,b-1,c;k,l-1) + c p(a+1,b,c-1;k,l-1) ]
\nonumber\\
& &\qquad\quad {}+ \rho c p(a+1,b+1,c-1;k,l),\nonumber
\end{eqnarray}
where $p(a,b,c;k,l) = 0$ if $a<0$, $b<0$, $c<0$, $k<0$, $l<0$,
$a=b=c=0$, or $k=l=0$. Equation (\ref{KL}) has a unique solution
satisfying the initial conditions
\[
p(1,0,0;k,l) = \delta_{k,1}\delta_{l,0},\qquad p(0,1,0;k,l) = \delta
_{k,0}\delta_{l,1}
\]
for $k,l = 0,1,\ldots, n$.

As with Golding's recursion, equation (\ref{KL}) can be solved
numerically, but quickly becomes computationally intractable with
growing $n$. The only exception is the special case of $\rho=\infty$,
for which the distribution is given by
the product of (\ref{eq:K_n}) for each locus.
In what follows, we use the following notation in writing
an asymptotic series for $p(a,b,c;k,l)$.
\begin{definition}
For loci $A$ and $B$, respectively, we define the analogues of (\ref
{eq:K_n}) as
%
%
\begin{equation}
\label{K}
p^A(a;k) = \frac{s(a,k)\theta_A^k}{(\theta_A)_{a}}
\end{equation}
and
%
%
\begin{equation}
\label{L}
p^B(b;l) = \frac{s(b,l)\theta_B^l}{(\theta_B)_{b}},
\end{equation}
where $s(a,k)$ and $s(b,l)$ are the Stirling numbers of the first kind.
\end{definition}

We pose the expansion
%
%
\begin{equation}
\label{expansionKL}
p(a,b,c;k,l) = p_0(a,b,c;k,l) + \frac{p_1(a,b,c;k,l)}{\rho} + O
\biggl(\frac{1}{\rho^2} \biggr)
\end{equation}
for large $\rho$. Then, in Section \ref{sec:proof_leadingKL},
we prove the following result for the
zeroth-order term.
%
%
\begin{proposition}
\label{prop:leadingKL}
For an asymptotic expansion of the form (\ref{expansionKL})
satisfying the recursion (\ref{KL}), $p_0(a,b,c;k,l)$ is given by
%
%
\begin{equation}
\label{leadingKL}
p_0(a,b,c;k,l) = p^A(a+c;k)p^B(b+c;l).
\end{equation}
\end{proposition}
%

Similar to Lemma \ref{lemma:first}, we obtain the following
vanishing result for the first-order term in the case of $c=0$.
%
%
\begin{lemma}
\label{lemma:firstKL}
For an asymptotic expansion of the form (\ref{expansionKL})
satisfying the recursion (\ref{KL}), we have
\[
p_1(a,b,0;k,l) = 0.
\]
\end{lemma}
%

Using this lemma, it is then possible to obtain the following
result for an arbitrary $c$.
%
\begin{proposition}
\label{prop:firstKL}
For an asymptotic expansion of the form (\ref{expansionKL})
satisfying the recursion (\ref{KL}), $p_1(a,b,c;k,l)$ is given by
%
%
\begin{eqnarray}
\label{firstorderKL}
p_1(a,b,c;k,l) &=& \frac{c(c-1)}{2} [ p^A(a+c;k) - p^A
(a+c-1;k) ]\nonumber\\[-8pt]\\[-8pt]
&&{}\times[ p^B(b+c;l) - p^B(b+c-1;l) ].\nonumber
\end{eqnarray}
\end{proposition}

Proofs of Proposition \ref{prop:firstKL} and Lemma \ref{lemma:firstKL}
are provided in Sections \ref{sec:proof_firstKL} and \ref
{sec:proof_lemma_firstKL}, respectively.

\section{Proofs of main results}
In what follows, we provide proofs of the results mentioned in the
previous two sections.

\subsection[Proof of Proposition 3.1]{Proof of Proposition \protect\ref{prop:leading}}
\label{sec:proof_leading}
First, assume $c > 0$. Substitute the expansion (\ref{expansion}) into
Golding's recursion (\ref{golding}),
divide by $\rho c$ and let $\rho\to\infty$. We are then left with
%
%
\begin{equation}
\label{unwrap0}
q_0({\mathbf{a}},{\mathbf{b}},{\mathbf{c}}) = \sum_{i=1}^K\sum
_{j=1}^L \frac{c_{ij}}{c}q_0({\mathbf{a}}+{\mathbf{e}}_i,{\mathbf
{b}}+{\mathbf{e}}
_j,{\mathbf{c}}-{\mathbf{e}}_{ij}).
\end{equation}
Now, applying (\ref{unwrap0}) repeatedly gives
\[
q_0({\mathbf{a}},{\mathbf{b}},{\mathbf{c}}) = \sum_{\mathrm{orderings}}
\frac{\prod_{(i,j)\in[K]\times[L]}c_{ij}!}{c!}
q_0({\mathbf{a}}+{\mathbf{c}}_A,{\mathbf{b}}+{\mathbf{c}}_B
,\mathbf{0}),
\]
where the summation is over all distinct orderings of the $c$ gametes
with multiplicity
${\mathbf{c}}=(c_{ij})$.
There are $\frac{c!}{\prod_{(i,j)}c_{ij}!}$ such orderings and since
the summand is independent of the ordering, we conclude
%
%
\begin{equation}
\label{unwrap0b}
q_0({\mathbf{a}},{\mathbf{b}},{\mathbf{c}}) = q_0({\mathbf
{a}}+{\mathbf{c}}_A,{\mathbf{b}}+{\mathbf{c}}_B,\mathbf{0}).
\end{equation}
Clearly, (\ref{unwrap0b}) also holds for $c = 0$.
From a
coalescent perspective, this equation tells us that any
gamete with specified alleles (i.e., ``carrying ancestral material'')
at both loci must undergo
recombination instantaneously backward in time.

Now, by substituting the asymptotic expansion (\ref{expansion}) with
${\mathbf{c}}= \mathbf{0}$ into Golding's recursion (\ref{golding})
and letting $\rho\to\infty$,
we obtain
%
%
\begin{eqnarray}
\label{goldingc=0preliminary}\quad
&&[n(n-1) + \theta_Aa + \theta_Bb ]q_0({\mathbf
{a}},{\mathbf{b}},\mathbf{0}) \nonumber\\
& &\qquad= \sum_{i=1}^K a_i(a_i-1)q_0({\mathbf{a}}-{\mathbf{e}}_i,{\mathbf
{b}},\mathbf{0}) + \sum_{j=1}^L b_j(b_j-1)q_0({\mathbf{a}},{\mathbf
{b}}-{\mathbf{e}}_j,\mathbf{0}) \nonumber\\[-8pt]\\[-8pt]
& &\qquad\quad {}+ 2\sum_{i=1}^K\sum_{j=1}^L a_i b_j q_0({\mathbf{a}}-{\mathbf{e}}
_i,{\mathbf{b}}-{\mathbf{e}}_j,{\mathbf{e}}_{ij}) \nonumber\\
& &\qquad\quad {}+ \theta_A\sum_{i=1}^K \delta_{a_i,1}q_0({\mathbf
{a}}-{\mathbf{e}}
_i,{\mathbf{b}},\mathbf{0}) + \theta_B\sum_{j=1}^L\delta
_{b_j,1}q_0({\mathbf{a}},{\mathbf{b}}-{\mathbf{e}}_j,\mathbf{0}).\nonumber
\end{eqnarray}
%
Equation (\ref{unwrap0b}) implies
$q_0({\mathbf{a}}-{\mathbf{e}}_i,{\mathbf{b}}-{\mathbf
{e}}_j,{\mathbf{e}}
_{ij})=q_0({\mathbf{a}},{\mathbf{b}},\mathbf{0})$, so
with a bit of rearranging we are left with
%
%
\begin{eqnarray}
\label{goldingc=0}\quad
&& [a(a+\theta_A-1) + b(b+\theta_B-1)
]q_0({\mathbf
{a}},{\mathbf{b}},\mathbf{0}) \nonumber\\
& &\qquad= \sum_{i=1}^K a_i(a_i-1)q_0({\mathbf{a}}-{\mathbf{e}}_i,{\mathbf
{b}},\mathbf{0}) + \sum_{j=1}^L b_j(b_j-1)q_0({\mathbf{a}},{\mathbf
{b}}-{\mathbf{e}}_j,\mathbf{0}) \\
& &\qquad\quad {}+ \theta_A\sum_{i=1}^K \delta_{a_i,1}q_0({\mathbf
{a}}-{\mathbf{e}}
_i,{\mathbf{b}},\mathbf{0}) + \theta_B\sum_{j=1}^L\delta
_{b_j,1}q_0({\mathbf{a}},{\mathbf{b}}-{\mathbf{e}}_j,\mathbf{0})\nonumber
\end{eqnarray}
%
with boundary conditions $q_0({\mathbf{e}}_i,\mathbf{0},\mathbf{0}) =
q_0(\mathbf{0},{\mathbf{e}}_j,\mathbf{0}) = 1$ for all $i\in[K]$
and \mbox{$j\in[L]$}.
Noting that (\ref{goldingc=0}) is the sum of two independent recursions
of the form~(\ref{one-locus-recursion}), one for each locus and each with
appropriate boundary condition, we conclude that $q_0({\mathbf
{a}},{\mathbf{b}},\mathbf{0})$ is
given by
%
%
\begin{equation}
\label{ESFx2}
q_0({\mathbf{a}},{\mathbf{b}},\mathbf{0}) = q^A({\mathbf{a}})q^B
({\mathbf{b}}),
\end{equation}
a product of two (ordered) ESFs. It is straightforward to verify that
(\ref{ESFx2}) satisfies~(\ref{goldingc=0}). Finally, using
(\ref{unwrap0b}) and (\ref {ESFx2}), we arrive at (\ref{leading}).

\subsection[Proof of Theorem 3.1]{Proof of Theorem \protect\ref{thm:first}}
\label{sec:proof_first}

First, assume $c >
0$. Substitute the asymptotic expansion (\ref{expansion}) into
Golding's recursion
(\ref{golding}), eliminate terms of order $\rho$ by
applying~(\ref{unwrap0}), and let $\rho\to\infty$. After applying (\ref{unwrap0b})
to the remaining terms and invoking~(\ref{goldingc=0}),
with some rearrangement we obtain
\begin{eqnarray*}
&&c q_1({\mathbf{a}},{\mathbf{b}},{\mathbf{c}}) - \sum
_{i=1}^K\sum_{j=1}^L c_{ij} q_1({\mathbf{a}}+{\mathbf{e}}_i,{\mathbf
{b}}+{\mathbf{e}}_j,{\mathbf{c}}-{\mathbf{e}}_{ij}) \nonumber
\\
& &\qquad= c(c-1)q_0({\mathbf{a}}+{\mathbf{c}}_A,{\mathbf{b}}+{\mathbf
{c}}_B,\mathbf{0})\\
&&\qquad\quad{}-\sum_{i=1}^K c_{i\cdot}(c_{i\cdot}-1)q_0({\mathbf{a}}+{\mathbf
{c}}_A-{\mathbf{e}}_i,{\mathbf{b}}+{\mathbf{c}}_B,\mathbf{0})
\nonumber
\\
& &\qquad\quad {}-\sum_{j=1}^L c_{\cdot j}(c_{\cdot j}-1)q_0({\mathbf
{a}}+{\mathbf{c}}_A,{\mathbf{b}}+{\mathbf{c}}_B-{\mathbf
{e}}_j,\mathbf
{0}) \nonumber\\
& &\qquad\quad {}+ \sum_{i=1}^K\sum_{j=1}^L
c_{ij}(c_{ij}-1)q_0({\mathbf{a}}+{\mathbf{c}}_A-{\mathbf
{e}}_i,{\mathbf
{b}}+{\mathbf{c}}_B-{\mathbf{e}}_j,\mathbf{0}).
\end{eqnarray*}
Now, by utilizing (\ref{ESFx2}), this can be written in the form
%
%
\begin{equation}
\label{q1+f}
q_1({\mathbf{a}},{\mathbf{b}},{\mathbf{c}}) = f({\mathbf
{a}},{\mathbf{b}},{\mathbf{c}}) + \sum_{i=1}^K\sum_{j=1}^L \frac
{c_{ij}}{c} q_1({\mathbf{a}}+{\mathbf{e}}_i,{\mathbf{b}}+{\mathbf
{e}}_j,{\mathbf
{c}}-{\mathbf{e}}_{ij}),
\end{equation}
where
%
%
\begin{eqnarray}
\label{eq:f}\qquad
f({\mathbf{a}},{\mathbf{b}},{\mathbf{c}}) &:=& (c-1)q^A({\mathbf
{a}}+{\mathbf{c}}_A)q^B({\mathbf{b}}+{\mathbf{c}}_B) \nonumber
\\
& & {}-q^B({\mathbf{b}}+{\mathbf{c}}_B)\sum_{i=1}^K \frac
{c_{i\cdot}(c_{i\cdot}-1)}{c}q^A({\mathbf{a}}+{\mathbf{c}}_A
-{\mathbf{e}}_i) \nonumber\\[-8pt]\\[-8pt]
& & {}-q^A({\mathbf{a}}+{\mathbf{c}}_A)\sum_{j=1}^L \frac
{c_{\cdot j}(c_{\cdot j}-1)}{c}q^B({\mathbf{b}}+{\mathbf{c}}_B
-{\mathbf{e}}_j) \nonumber\\
& & {}+ \sum_{i=1}^K\sum_{j=1}^L \frac{c_{ij}(c_{ij}-1)}{c}q^A
({\mathbf{a}}+{\mathbf{c}}_A-{\mathbf{e}}_i)q^B({\mathbf
{b}}+{\mathbf
{c}}_B-{\mathbf{e}}_j).\nonumber
\end{eqnarray}
Above, we assumed $c > 0$.
We define $f({\mathbf{a}},{\mathbf{b}},{\mathbf{c}})=0$ if ${\mathbf
{c}}=\mathbf{0}$.
Iterating the recursion~(\ref{q1+f}), we may write
$q_1({\mathbf{a}},{\mathbf{b}},{\mathbf{c}})$ as
%
\begin{eqnarray*}
q_1({\mathbf{a}},{\mathbf{b}},{\mathbf{c}}) &=& f({\mathbf
{a}},{\mathbf{b}},{\mathbf{c}})\\
&&{} + \sum_{i=1}^K \sum_{j=1}^L
\frac{c_{ij}}{c} \Biggl[f({\mathbf{a}}+{\mathbf{e}}_i,{\mathbf
{b}}+{\mathbf{e}}
_j,{\mathbf{c}}-{\mathbf{e}}_{ij})
\\
& &\hspace*{62.4pt} {}+ \sum_{i'=1}^K \sum_{j'=1}^L \frac{c_{i'j'} -
\delta_{ii'}\delta_{jj'}}{c-1}\\
&&\hspace*{110.7pt}{}\times q_1({\mathbf{a}}+{\mathbf
{e}}_i+{\mathbf{e}}_{i'},{\mathbf
{b}}+{\mathbf{e}}_j+{\mathbf{e}}_{j'},\\
&&\hspace*{190pt}{\mathbf{c}}-{\mathbf
{e}}_{ij}-{\mathbf{e}}_{i'j'})\Biggr].
\end{eqnarray*}
%
Similarly, repeatedly iterating (\ref{q1+f}) yields
%
%
\begin{eqnarray}
\label{expandf}\qquad
q_1({\mathbf{a}},{\mathbf{b}},{\mathbf{c}})
&=&
q_1({\mathbf{a}}+{\mathbf{c}}_A,{\mathbf{b}}+{\mathbf{c}}_B
,\mathbf{0}) + f({\mathbf{a}},{\mathbf{b}},{\mathbf{c}}) \nonumber\\
& & {}+ \sum_{i_1 j_1}\frac{c_{i_1j_1}}{c}f({\mathbf{a}}+{\mathbf{e}}
_{i_1},{\mathbf{b}}+{\mathbf{e}}_{j_1},{\mathbf{c}}-{\mathbf{e}}_{i_1j_1})
\nonumber\\
& & {}+ \sum_{i_1 j_1,i_2 j_2}\frac{c_{i_1 j_1}}{c}\frac{c_{i_2 j_2}
- \delta_{{i_1 j_1},{i_2 j_2}}}{c-1} \\
& &\hspace*{42.5pt}{} \times f({\mathbf
{a}}+{\mathbf{e}}_{i_1}+{\mathbf{e}}_{i_2},{\mathbf{b}}+{\mathbf
{e}}_{j_1}+{\mathbf{e}}
_{j_2},{\mathbf{c}}-{\mathbf{e}}_{i_1 j_1}-{\mathbf{e}}_{i_2 j_2})
\nonumber\\
& & {}+ \cdots+ \sum_{i_1 j_1,\ldots,i_c j_c}\frac{\prod
_{ij}c_{ij}!}{c!}f({\mathbf{a}}+{\mathbf{c}}_A,{\mathbf
{b}}+{\mathbf{c}}_B,\mathbf{0}).\nonumber
\end{eqnarray}
The key observation is that the right-hand side of (\ref{expandf}) has
a nice
probabilistic interpretation which allows us to obtain a closed-form formula.
To be more precise, consider the first summation
\[
\sum_{i_1 j_1}\frac{c_{i_1 j_1}}{c}f({\mathbf{a}}+{\mathbf{e}}
_{i_1},{\mathbf{b}}+{\mathbf{e}}_{j_1},{\mathbf{c}}-{\mathbf
{e}}_{i_1 j_1}).
\]
For a fixed sample configuration ${\mathbf{c}}$, this can be
interpreted as the sum over all
possible ways of throwing away a gamete at random and calculating $f$
based on the remaining subsample, which we will denote
${\mathbf{c}}^{(c-1)}$. Equivalently, it is the expected value of $f$ with
respect to subsampling without replacement $c-1$ of the gametes in
${\mathbf{c}}$. Write this as
\[
\mathbb{E}\bigl[f\bigl(\mathbf{A}^{(c-1)},\mathbf{B}^{(c-1)},\mathbf{C}^{(c-1)}\bigr)\bigr],
\]
where $\mathbf{C}^{(c-1)}$ is the random subsample obtained by sampling
without replacement $c-1$ gametes from ${\mathbf{c}}$, and $\mathbf{A}
^{(c-1)} := {\mathbf{a}}
+ {\mathbf{c}}_A- \mathbf{C}_A^{(c-1)}$, $\mathbf{B}^{(c-1)} :=
{\mathbf
{b}}+ {\mathbf{c}}_B-
\mathbf{C}_B^{(c-1)}$. Note that once the subsample ${\mathbf
{c}}^{(c-1)}$ is
obtained, then ${\mathbf{a}}^{(c-1)}$ and ${\mathbf{b}}^{(c-1)}$ are fully
specified. More generally, consider the $(c-m)$th sum in
(\ref{expandf}). A particular
term in the summation corresponds to an ordering of $c-m$ gametes in
${\mathbf{c}}$, which, when removed leave a subsample ${\mathbf
{c}}^{(m)}$. With
respect to this subsample, the summand is
\[
\prod_{i=1}^K \prod_{j=1}^L \frac{c_{ij}!}{c_{ij}^{(m)}!}\frac
{m!}{c!}f\bigl({\mathbf{a}}^{(m)},{\mathbf{b}}^{(m)},{\mathbf{c}}^{(m)}\bigr)
\]
and for each such subsample ${\mathbf{c}}^{(m)}$ there are $
{c-m\choose{\mathbf{c}}-{\mathbf{c}}^{(m)}}$ distinct
orderings of the remaining types in ${\mathbf{c}}$, with each ordering
contributing the same amount to the
sum.
Here, ${c-m\choose{\mathbf{c}}-{\mathbf{c}}^{(m)}}$ denotes the multinomial
coefficient:
\[
\pmatrix{c-m\cr{\mathbf{c}}-{\mathbf{c}}^{(m)}} = \frac
{(c-m)!}{\prod
_{i=1}^K \prod_{j=1}^L (c_{ij}-c_{ij}^{(m)})!} .
\]
Gathering identical terms, the $(c-m)$th sum in (\ref{expandf}) can
therefore be written
over all distinct subsamples of ${\mathbf{c}}$ of size $m$:
\begin{eqnarray*}
&&\sum_{{\mathbf{c}}^{(m)}} \pmatrix{c-m\cr{\mathbf
{c}}-{\mathbf{c}}^{(m)}} \prod_{i=1}^K \prod_{j=1}^L\frac
{c_{ij}!}{c_{ij}^{(m)}!}\frac{m!}{c!}f\bigl({\mathbf{a}}^{(m)},{\mathbf
{b}}^{(m)},{\mathbf{c}}^{(m)}\bigr)\\
&&\qquad= \sum_{{\mathbf
{c}}^{(m)}}\frac{1}{{c\choose m}} \prod_{i=1}^K \prod_{j=1}^L
\pmatrix{c_{ij}\cr c_{ij}^{(m)}} f\bigl({\mathbf{a}}^{(m)},{\mathbf
{b}}^{(m)},{\mathbf{c}}^{(m)}\bigr) \\
&&\qquad= \mathbb{E}\bigl[f\bigl(\mathbf{A}^{(m)},\mathbf{B}^{(m)},\mathbf{C}^{(m)}\bigr)\bigr],
\end{eqnarray*}
where, for a fixed $m$, $\mathbf{C}^{(m)}=(C_{ij}^{(m)})$ is a multivariate
hypergeometric$(c,{\mathbf{c}},m)$ random variable; that is,
\[
\mathbb{P} \biggl(\bigcap_{(i,j)\in[K]\times[L]} \bigl[C_{ij}^{(m)} =
c_{ij}^{(m)} \bigr] \biggr) =
\frac{1}{{c\choose m} }
\prod_{(i,j)\in[K]\times[L]} \pmatrix{c_{ij}\cr c_{ij}^{(m)}}.
\]
Furthermore, marginally we
have
\begin{eqnarray*}
C_{ij}^{(m)} &\sim& \operatorname{hypergeometric}(c,c_{ij},m), \\
C_{i\cdot}^{(m)} &\sim& \operatorname{hypergeometric}(c,c_{i\cdot
},m), \\
C_{\cdot j}^{(m)} &\sim& \operatorname{hypergeometric}(c,c_{\cdot j},m).
\end{eqnarray*}
In summary, (\ref{expandf}) can be written as
%
%
\begin{equation}
\label{q1concise}
q_1({\mathbf{a}},{\mathbf{b}},{\mathbf{c}}) = q_1({\mathbf
{a}}+{\mathbf{c}}_A,{\mathbf{b}}+{\mathbf{c}}_B,\mathbf{0}) +
\sum_{m=1}^c \mathbb{E}\bigl[f\bigl(\mathbf{A}^{(m)},\mathbf{B}^{(m)},\mathbf
{C}^{(m)}\bigr)\bigr].
\end{equation}
According to Lemma \ref{lemma:first}, the first term
$q_1({\mathbf{a}}+{\mathbf{c}}_A,{\mathbf{b}}+{\mathbf{c}}_B
,\mathbf{0})$ vanishes, so we are left with
%
%
\begin{equation}
\label{subsampling}
q_1({\mathbf{a}},{\mathbf{b}},{\mathbf{c}}) = \sum_{m=1}^c \mathbb{E}
\bigl[f\bigl(\mathbf{A}^{(m)},\mathbf{B}^{(m)},\mathbf{C}^{(m)}\bigr)\bigr].
\end{equation}

Finally, since $\mathbf{A}^{(m)}+\mathbf{C}_A^{(m)} = {\mathbf
{a}}+{\mathbf
{c}}_A$ and
$\mathbf{B}^{(m)}+\mathbf{C}_B^{(m)} = {\mathbf{b}}+{\mathbf{c}}_B$,
(\ref{eq:f}) and (\ref{subsampling}) together imply
%
\begin{eqnarray*}
&&q_1({\mathbf{a}},{\mathbf{b}},{\mathbf{c}}) \\
&&\qquad=\sum_{m=1}^c \Biggl[(m-1)q^A({\mathbf{a}}+{\mathbf{c}}_A)q^B
({\mathbf{b}}+{\mathbf{c}}_B)
\\
& &\qquad\quad\hspace*{21.2pt}{}
-q^B({\mathbf{b}}+{\mathbf{c}}_B)\frac{1}{m}\sum_{i=1}^K \mathbb{E}
\bigl[C_{i\cdot}^{(m)}\bigl(C_{i\cdot}^{(m)}-1\bigr)\bigr]q^A({\mathbf{a}}+{\mathbf
{c}}_A-{\mathbf{e}}_i) \\
& &\qquad\quad\hspace*{21.2pt}{} - q^A({\mathbf{a}}+{\mathbf{c}}_A)\frac{1}{m}\sum_{j=1}^L
\mathbb{E}\bigl[C_{\cdot j}^{(m)}\bigl(C_{\cdot j}^{(m)}-1\bigr)\bigr]q^B({\mathbf
{b}}+{\mathbf{c}}_B-{\mathbf{e}}_j) \\
& &\qquad\quad\hspace*{21.2pt}{} + \frac{1}{m}\sum_{i=1}^K \sum_{j=1}^L \mathbb{E}
\bigl[C_{ij}^{(m)}\bigl(C_{ij}^{(m)}-1\bigr)\bigr]q^A({\mathbf{a}}+{\mathbf{c}}_A
-{\mathbf{e}}_i)q^B({\mathbf{b}}+{\mathbf{c}}_B-{\mathbf{e}}_j)\Biggr].
\end{eqnarray*}
The moments in this equation are easy to compute and one can sum them
over $m$ to
obtain the desired result (\ref{firstorder}).

\subsection[Proof of Lemma 3.1]{Proof of Lemma \protect\ref{lemma:first}}
\label{sec:proof_lemma_first}

First, note that for any sample $({\mathbf{a}},{\mathbf{b}},{\mathbf
{c}})$ and any subsample
of the form $({\mathbf{a}}^{(1)},{\mathbf{b}}^{(1)},{\mathbf
{c}}^{(1)})$, we have
$f({\mathbf{a}}^{(1)},{\mathbf{b}}^{(1)},{\mathbf{c}}^{(1)}) = 0$,
since every term on
right-hand side of (\ref{eq:f}) has a vanishing coefficient.
So, equation
(\ref{q1concise}) implies
%
%
\begin{equation}
\label{unwrap1b}
q_1({\mathbf{a}}-{\mathbf{e}}_i,{\mathbf{b}}-{\mathbf
{e}}_j,{\mathbf{e}}_{ij}) =
q_1({\mathbf{a}},{\mathbf{b}},\mathbf{0})
\end{equation}
for any $(i,j)\in[K]\times[L]$. Now, substitute the asymptotic expansion
(\ref{expansion}) with ${\mathbf{c}}= \mathbf{0}$ into Golding's recursion
(\ref{golding}).
Note that terms of order $\rho$ are absent since ${\mathbf{c}}=
\mathbf{0}$.
Eliminate terms with coefficients independent of $\rho$ by applying
(\ref{goldingc=0preliminary}), multiply both sides of the recursion by
$\rho$,
and let
$\rho\to\infty$ to obtain the following:
\begin{eqnarray*}
&&[n(n-1) + \theta_Aa + \theta_Bb]q_1({\mathbf{a}},{\mathbf
{b}},\mathbf{0}) \\
& &\qquad= \sum_{i=1}^K a_i(a_i-1)q_1({\mathbf{a}}-{\mathbf{e}}_i,{\mathbf
{b}},\mathbf{0}) + \sum_{j=1}^L b_j(b_j-1)q_1({\mathbf{a}},{\mathbf
{b}}-{\mathbf{e}}_j,\mathbf{0}) \\
& &\qquad\quad {}+ 2\sum_{i=1}^K \sum_{j=1}^L a_i b_j q_1({\mathbf{a}}-{\mathbf{e}}
_i,{\mathbf{b}}-{\mathbf{e}}_j,{\mathbf{e}}_{ij}) \\
& &\qquad\quad {}+ \theta_A\sum_{i=1}^K \delta_{a_i,1}q_1({\mathbf
{a}}-{\mathbf{e}}
_i,{\mathbf{b}},\mathbf{0}) + \theta_B\sum_{j=1}^L\delta
_{b_j,1}q_1({\mathbf{a}},{\mathbf{b}}-{\mathbf{e}}_j,\mathbf{0})
\end{eqnarray*}
with boundary conditions $q_1({\mathbf{e}}_i,\mathbf{0},\mathbf{0}) =
q_1(\mathbf{0},{\mathbf{e}}_j,\mathbf{0}) = 0$ for all $i\in[K]$
and $j\in[L]$.
This equation can be made
strictly recursive by applying (\ref{unwrap1b}) to
$q_1({\mathbf{a}}-{\mathbf{e}}_i,{\mathbf{b}}-{\mathbf
{e}}_j,{\mathbf{e}}_{ij})$. It
therefore follows
from the boundary conditions (for example, by induction) that
$q_1({\mathbf{a}},{\mathbf{b}},\mathbf{0}) = 0$.

\subsection[Proof of Theorem 3.2]{Proof of Theorem \protect\ref{thm:second}}
\label{sec:proof_second}

Here, we provide only an outline of a proof; details are
similar to the proof of Theorem \ref{thm:first}.
Substitute the asymptotic expansion
(\ref{expansion}) into Golding's recursion (\ref{golding}),
eliminate terms
with coefficients proportional to $\rho$ or independent of $\rho$.
Then, multiply both sides of the recursion by $\rho$
and let $\rho\to\infty$ to obtain
%
%
\begin{eqnarray}
\label{secondrec}
&&c q_2({\mathbf{a}},{\mathbf{b}},{\mathbf{c}}) - \sum
_{i=1}^K \sum_{j=1}^L c_{ij}q_2({\mathbf{a}}+{\mathbf{e}}_i,{\mathbf
{b}}+{\mathbf{e}}_j,{\mathbf{c}}-{\mathbf{e}}_{ij}) \nonumber\\
& &\qquad= \sum_{i=1}^K a_i(a_i-1+2c_{i\cdot})q_1({\mathbf
{a}}-{\mathbf{e}}_i,{\mathbf{b}},{\mathbf{c}})\nonumber\\
&&\qquad\quad{} + \sum_{j=1}^L
b_j(b_j-1+2c_{\cdot j})q_1({\mathbf{a}},{\mathbf{b}}-{\mathbf
{e}}_j,{\mathbf
{c}}) \nonumber\\
& &\qquad\quad{} + \sum_{i=1}^K \sum_{j=1}^L
[ c_{ij}(c_{ij}-1)q_1({\mathbf{a}},{\mathbf{b}},{\mathbf
{c}}-{\mathbf{e}}_{ij})\nonumber\\
&&\hspace*{78.7pt}{} + 2 a_i b_j q_1({\mathbf{a}}-{\mathbf
{e}}_i,{\mathbf
{b}}-{\mathbf{e}}_j,{\mathbf{c}}+{\mathbf{e}}_{ij}) ]
\nonumber\\[-8pt]\\[-8pt]
& &\qquad\quad{} + \theta_A\sum_{i=1}^K \Biggl[\sum_{j=1}^L \delta
_{a_i+c_{i\cdot},1}\delta_{c_{ij},1}q_1({\mathbf{a}},{\mathbf
{b}}+{\mathbf{e}}
_j,{\mathbf{c}}-{\mathbf{e}}_{ij})\nonumber\\
&&\qquad\quad\hspace*{93.3pt}{} + \delta_{a_i,1}\delta
_{c_{i\cdot
},0}q_1({\mathbf{a}}-{\mathbf{e}}_i,{\mathbf{b}},{\mathbf
{c}}) \Biggr]
\nonumber\\
& &\qquad\quad{} + \theta_B\sum_{j=1}^L \Biggl[\sum_{i=1}^K \delta
_{b_j+c_{\cdot j},1}\delta_{c_{ij},1}q_1({\mathbf{a}}+{\mathbf
{e}}_i,{\mathbf
{b}},{\mathbf{c}}-{\mathbf{e}}_{ij})\nonumber\\
&&\qquad\quad\hspace*{89.3pt}{} + \delta_{b_j,1}\delta
_{c_{\cdot
j},0}q_1({\mathbf{a}},{\mathbf{b}}-{\mathbf{e}}_j,{\mathbf
{c}}) \Biggr]
\nonumber\\
& &\qquad\quad{} - [n(n-1) + \theta_A(a+c) + \theta
_B(b+c)]q_1({\mathbf
{a}},{\mathbf{b}},{\mathbf{c}}). \nonumber
\end{eqnarray}
By substituting our expression (\ref{firstorder}) for
$q_1({\mathbf{a}},{\mathbf{b}},{\mathbf{c}})$, the right-hand side can
be expressed as a function $g({\mathbf{a}},{\mathbf{b}},{\mathbf{c}})$
which is completely known but rather cumbersome to write down. As in
the proof of Theorem \ref{thm:first}, the same ``unwrapping'' maneuver
can be applied to rearrange (\ref{secondrec}) into the form
(\ref{q2concise}), where
\[
\sigma({\mathbf{a}},{\mathbf{b}},{\mathbf{c}}) = \sum_{m=1}^c
\mathbb{E}\bigl[g\bigl(\mathbf{A}^{(m)},\mathbf{B}^{(m)},\mathbf{C}^{(m)}\bigr)\bigr].
\]
This time
$\mathbb{E}[g(\mathbf{A}^{(m)},\mathbf{B}^{(m)},\mathbf{C}^{(m)})]$ is
a function of fourth-order moments of the multivariate hypergeometric
distribution. The formula shown in the \hyperref[appendixA]{Appendix}
is obtained by evaluating the expectations and summing over $m$.

We now show that $q_2({\mathbf{a}},{\mathbf{b}},\mathbf{0})$ satisfies
the recursion shown in (\ref{q2(a,b,0)}). We will use the fact that for
a sample
$({\mathbf{a}}-{\mathbf{e}}_i,{\mathbf{b}}-{\mathbf{e}}_j,{\mathbf
{e}}_{ij})$, we have
%
%
\begin{eqnarray}
\label{g(a1,b1,c1)}
&& g({\mathbf{a}}-{\mathbf{e}}_i,{\mathbf{b}}-{\mathbf{e}}_j,{\mathbf
{e}}_{ij}) \nonumber\\
&&\qquad=
2(a-1)(b-1)q^A({\mathbf{a}})q^B({\mathbf{b}}) \nonumber\\
&&\qquad\quad {}- 2(b-1)(a_i-1)q^A({\mathbf{a}}-{\mathbf{e}}_i)q^B({\mathbf{b}})
\\
&&\qquad\quad {}- 2(a-1)(b_j-1)q^A({\mathbf{a}})q^B({\mathbf{b}}-{\mathbf{e}}_j)
\nonumber\\
&&\qquad\quad {}+ 2(a_i-1)(b_j-1)q^A({\mathbf{a}}-{\mathbf{e}}_i)q^B({\mathbf
{b}}-{\mathbf{e}}_j).\nonumber
\end{eqnarray}
For an arbitrary ${\mathbf{c}}$, $g({\mathbf{a}},{\mathbf
{b}},{\mathbf{c}})$ is much more complicated.

Now, one can adopt the approach used in the proof of Lemma \ref
{lemma:first} to obtain a strict recursion for
$q_2({\mathbf{a}}+{\mathbf{c}}_A,{\mathbf{b}}+{\mathbf{c}}_B
,\mathbf{0})$. First, note that (\ref{q2concise}) and
(\ref{g(a1,b1,c1)}) imply
%
%
\begin{eqnarray}
\label{unwrap2}
&&q_2({\mathbf{a}}-{\mathbf{e}}_i,{\mathbf{b}}-{\mathbf
{e}}_j,{\mathbf{e}}_{ij}) \nonumber\\
&&\qquad=
q_2({\mathbf{a}},{\mathbf{b}},\mathbf{0}) + \mathbb{E}\bigl[g\bigl((\mathbf
{A}-{\mathbf{e}}
_i)^{(1)},(\mathbf{B}-{\mathbf{e}}_j)^{(1)},{\mathbf
{e}}_{ij}^{(1)}\bigr)\bigr] \nonumber\\
&&\qquad= q_2({\mathbf{a}},{\mathbf{b}},\mathbf{0}) + g({\mathbf
{a}}-{\mathbf{e}}
_i,{\mathbf{b}}-{\mathbf{e}}_j,{\mathbf{e}}_{ij}) \nonumber\\
&&\qquad= q_2({\mathbf{a}},{\mathbf{b}},\mathbf{0}) + 2(a-1)(b-1)q^A
({\mathbf{a}})q^B({\mathbf{b}}) \\\
&&\qquad\quad {}- 2(b-1)(a_i-1)q^A({\mathbf{a}}-{\mathbf{e}}_i)q^B({\mathbf{b}})
\nonumber\\
&&\qquad\quad {}- 2(a-1)(b_j-1)q^A({\mathbf{a}})q^B({\mathbf{b}}-{\mathbf{e}}_j)
\nonumber\\
&&\qquad\quad {}+ 2(a_i-1)(b_j-1)q^A({\mathbf{a}}-{\mathbf{e}}_i)q^B({\mathbf
{b}}-{\mathbf{e}}_j).\nonumber
\end{eqnarray}
As before, substitute the asymptotic expansion
(\ref{expansion})
for ${\mathbf{c}}=\mathbf{0}$ into Golding's recursion (\ref{golding}),
eliminate terms with coefficients independent of
$\rho$ or proportional to $\rho^{-1}$, and let $\rho\to\infty$
to obtain
\begin{eqnarray*}
&&[n(n-1) + \theta_Aa + \theta_Bb]q_2({\mathbf{a}},{\mathbf
{b}},\mathbf{0}) \\
&&\qquad= \sum_{i=1}^K a_i(a_i-1)q_2({\mathbf{a}}-{\mathbf{e}}_i,{\mathbf
{b}},\mathbf{0}) + \sum_{j=1}^L b_j(b_j-1)q_2({\mathbf{a}},{\mathbf
{b}}-{\mathbf{e}}_j,\mathbf{0}) \\
&&\qquad\quad {}+ 2\sum_{i=1}^K \sum_{j=1}^L a_i b_j q_2({\mathbf{a}}-{\mathbf{e}}
_i,{\mathbf{b}}-{\mathbf{e}}_j,{\mathbf{e}}_{ij}) \\
&&\qquad\quad{} + \theta_A\sum_{i=1}^K \delta_{a_i,1}q_2({\mathbf
{a}}-{\mathbf{e}}
_i,{\mathbf{b}},\mathbf{0}) + \theta_B\sum_{j=1}^L\delta
_{b_j,1}q_2({\mathbf{a}},{\mathbf{b}}-{\mathbf{e}}_j,\mathbf{0})
\end{eqnarray*}
with boundary conditions $q_2({\mathbf{e}}_i,\mathbf{0},\mathbf{0}) =
q_2(\mathbf{0},{\mathbf{e}}_j,\mathbf{0}) = 0$ for all $i\in[K]$
and $j\in[L]$.
This equation can be made strictly recursive by applying
(\ref{unwrap2}) to $q_2({\mathbf{a}}-{\mathbf{e}}_i,{\mathbf
{b}}-{\mathbf{e}}
_j,{\mathbf{e}}_{ij})$. After
some simplification, this leads to the recursion
(\ref{q2(a,b,0)}).

\subsection[Proof of Proposition 4.1]{Proof of Proposition \protect\ref{prop:leadingKL}}
\label{sec:proof_leadingKL}
The proof is similar to the proof of Proposition~\ref{prop:leading},
working with the system (\ref{KL}) rather than Golding's recursion
(\ref{golding}). First, assume $c > 0$. Substitute the expansion (\ref
{expansionKL}) into the recursion (\ref{KL}), divide by $\rho c$ and
let $\rho\to \infty$. We are left with
\[
p_0(a,b,c;k,l) = p_0(a+1,b+1,c-1;k,l),
\]
which implies
%
%
\begin{equation}
\label{unwrapKL}
p_0(a,b,c;k,l) = p_0(a+c,b+c,0;k,l).
\end{equation}
Clearly, (\ref{unwrapKL}) also holds for $c = 0$.

Now, by substituting the asymptotic expansion (\ref{expansionKL}) with
$c=0$ into (\ref{KL}) and letting $\rho\to\infty$, we obtain
%
%
\begin{eqnarray}
\label{KLc=0preliminary}\quad
&& [n(n-1) + \theta_Aa + \theta_Bb ]p_0(a,b,0;k,l)
\nonumber\\
&&\qquad= a(a-1)p_0(a-1,b,0;k,l) + b(b-1)p_0(a,b-1,c;k,l)
\nonumber\\
&&\qquad\quad {}+ 2ab p_0(a-1,b-1,1;k,l) \\
&&\qquad\quad {}+ \theta_Aa p_0(a-1,b,0;k-1,l)\nonumber\\
&&\qquad\quad{} + \theta_Bb
p_0(a,b-1,0;k,l-1).\nonumber
\end{eqnarray}
After invoking (\ref{unwrapKL}) on $p_0(a-1,b-1,1;k,l)$ and
rearranging, we are left with
%
%
\begin{eqnarray}
\label{KLc=0}\quad
&& [a(a+\theta_A-1) + b(b+\theta_B-1)
]p_0(a,b,0;k,l)
\nonumber\\
&&\qquad = a(a-1)p_0(a-1,b,0;k,l) + b(b-1)p_0(a,b-1,0;k,l)
\nonumber\\[-8pt]\\[-8pt]
&&\qquad\quad{}+ \theta_Aa p_0(a-1,b,0;k-1,l)\nonumber\\
&&\qquad\quad{} + \theta_Bb
p_0(a,b-1,0;k,l-1)\nonumber
\end{eqnarray}
with boundary conditions $p_0(1,0,0;k,l) = \delta_{k,1}\delta
_{l,0}$ and
$p_0(0,1,0;k,l) = \delta_{k,0}\delta_{l,1}$. Equation (\ref
{KLc=0}) can be
expressed as a linear sum of two independent recursions:
\begin{eqnarray*}
(a-1+\theta_A)p_0^A(a;k) &=& (a-1)p_0^A(a-1;k) + \theta_Ap_0^A(a-1;k-1),
\\
(b-1+\theta_B)p_0^B(b;l) &=& (b-1)p_0^B(b-1;l) + \theta_Bp_0^B(b-1;l-1)
\end{eqnarray*}
with respective boundary conditions $p_0^A(1;k) = \delta_{k,1}$ and
$p_0^B(1;l) = \delta_{l,1}$. These recursions are precisely those
considered by Ewens [(\citeyear{ewe1972}), (21)],
with respective solutions
(\ref{K}) and (\ref{L}). Hence, $p_0^A(a;k) = p^A(a;k)$ and
$p_0^B(b;l) = p^B(b;l)$, and it is straightforward to verify that
$p^A(a;k)p^B(b;l)$ satisfies (\ref{KLc=0}). Substituting this
solution into (\ref{unwrapKL}), we arrive at (\ref{leadingKL}), as
required.

\subsection[Proof of Proposition 4.2]{Proof of Proposition \protect\ref{prop:firstKL}}
\label{sec:proof_firstKL}

First, assume
$c > 0$. Substitute the asymptotic expansion (\ref{expansionKL}) into
the recursion (\ref{KL}), eliminate terms with coefficients linear in
$\rho$ by applying
(\ref{unwrapKL}), and let $\rho\to\infty$. After applying
(\ref{unwrapKL}) to the remaining terms and invoking (\ref{KLc=0}),
with some rearrangement we obtain
%
%
\begin{eqnarray}
\label{firstrecKL}\qquad
&&p_1(a,b,c;k,l) - p_1(a+1,b+1,c-1;k,l) \nonumber\\
&&\qquad= (c-1)[p_0(a+c,b+c,0;k,l) - p_0(a+c-1,b+c,0;k,l)
\nonumber\\[-8pt]\\[-8pt]
&&\qquad\quad\hspace*{33pt}{}- p_0(a+c,b+c-1,0;k,l)\nonumber\\
&&\qquad\quad\hspace*{120.04pt}{} + p_0(a+c-1,b+c-1,0;k,l)].\nonumber
\end{eqnarray}
Applying the recursion repeatedly, this becomes
%
%
\begin{eqnarray}
\label{p1+f}
&&p_1(a,b,c;k,l) \nonumber\\
&&\qquad= p_1(a+c,b+c,0;k,l)\nonumber\\
&&\qquad\quad{} + [p_0(a+c,b+c,0;k,l) -
p_0(a+c-1,b+c,0;k,l)\nonumber\\[-8pt]\\[-8pt]
&&\qquad\quad\hspace*{14.05pt}{} - p_0(a+c,b+c-1,0;k,l) \nonumber\\
&&\qquad\quad\hspace*{100.5pt}{}+ p_0(a+c-1,b+c-1,0;k,l)]\nonumber\\
&&\qquad\quad\hspace*{10.6pt}{}\times\sum_{m=0}^{c-1} (c-1-m).\nonumber
\end{eqnarray}

According to Lemma \ref{lemma:firstKL}, the first term
$p_1(a+c,b+c,0;k,l)$ vanishes. Hence, since $p_0(a,b,c;k,l)$ is given
by (\ref{leadingKL}), the right-hand side of (\ref{p1+f}) is fully
known. With some rearrangement, we are left with (\ref{firstorderKL}).

\subsection[Proof of Lemma 4.1]{Proof of Lemma \protect\ref{lemma:firstKL}}
\label{sec:proof_lemma_firstKL}
First, note that (\ref{firstrecKL}) implies
%
%
\begin{equation}
\label{unwrap1KL}
p_1(a-1,b-1,1;k,l) = p_1(a,b,0;k,l).
\end{equation}
Now, substitute the asymptotic expansion (\ref{expansionKL}) with $c =
0$ into (\ref{KL}), eliminate leading-order terms by applying
(\ref{KLc=0preliminary}), and let $\rho\to\infty$. The result is
made strictly
recursive by invoking (\ref{unwrap1KL}), and we obtain
\begin{eqnarray*}
&&[a(a+\theta_A-1) + b(b+\theta_B-1)]p_1(a,b,0;k,l)
\\
&&\qquad= a(a-1)p_1(a-1,b,0;k,l)\\
&&\qquad\quad{} + b(b-1)p_1(a,b-1,0;k,l) \\
&&\qquad\quad{}+ \theta_Aa p_1(a-1,b,0;k-1,l)\\
&&\qquad\quad{} + \theta_Bb p_1(a,b-1,0;k,l-1)
\end{eqnarray*}
with boundary conditions $p_1(1,0,0;k,l) = p_1(0,1,0;k,l) = 0$, for
$k,l = 0,\break\ldots,n$. It therefore follows (e.g., by induction)
that $p_1(a,b,0;k,l) = 0$.

\begin{appendix}\label{appendixA}
\section*{Appendix: Expression for $\sigma({\mathbf{a}},{\mathbf
{b}},{\mathbf{c}})$}

We use $Q^A$ to denote $q^A({\mathbf{a}}+{\mathbf{c}}_A)$,
$Q_i^A$ to denote
$q^A({\mathbf{a}}+{\mathbf{c}}_A-{\mathbf{e}}_i)$,
$Q_{ik}^A$ to denote $q^A({\mathbf{a}}+{\mathbf{c}}_A-{\mathbf{e}}
_i-{\mathbf{e}}_k)$, and so on. Then
\begin{eqnarray*}
\label{sumEg}
&&\sigma({\mathbf{a}},{\mathbf{b}},{\mathbf{c}})
\nonumber\\
&&\qquad=\frac{c}{3} \biggl[\frac{(c-1)(c+1)(3c-2)}{8} + (c-1)(3a+3b+2c-1) +
6ab \biggr]Q^AQ^B\nonumber\\
&&\qquad\quad {}- \frac{\theta_A(c-1)}{2}\sum_{i=1}^K \delta_{a_i,0}\delta
_{c_{i\cdot
},1} Q_i^AQ^B- \frac{\theta_B(c-1)}{2}\sum_{j=1}^L \delta
_{b_j,0}\delta
_{c_{\cdot j},1} Q^AQ_j^B\nonumber\\
&&\qquad\quad {}+ \sum_{i=1}^K \biggl[\frac{\theta_A- c(c-3) + 2a + 4b - 4}{4}
c_{i\cdot}(c_{i\cdot}-1) \nonumber\\
&&\qquad\quad\hspace*{72pt}{} -(2b+c-1) c_{i\cdot}(a_i + c_{i\cdot
}-1) \biggr]Q_i^AQ^B\nonumber\\
&&\qquad\quad {}+ \frac{1}{2}\sum_{i=1}^K \biggl[\frac{\theta_A}{2}\delta
_{c_{i\cdot
},2} + \frac{5-6a_i-4c_{i\cdot}}{6} \biggr]c_{i\cdot}(c_{i\cdot
}-1)Q_{ii}^AQ^B\nonumber\\
&&\qquad\quad {}+ \sum_{j=1}^L \biggl[\frac{\theta_B- c(c-3) + 2b + 4a - 4}{4}
c_{\cdot j}(c_{\cdot j}-1) \nonumber\\
&&\qquad\quad\hspace*{71.83pt}{} -(2a+c-1) c_{\cdot j}(b_j + c_{\cdot
j}-1) \biggr]Q^AQ_j^B\nonumber\\
&&\qquad\quad {}+ \frac{1}{2}\sum_{j=1}^L \biggl[\frac{\theta_B}{2}\delta
_{c_{\cdot
j},2} + \frac{5-6b_j-4c_{\cdot j}}{6} \biggr]c_{\cdot j}(c_{\cdot
j}-1) Q^AQ_{jj}^B\nonumber\\
&&\qquad\quad
{}+\sum_{i,k=1}^K \frac{c_{i\cdot}(c_{i\cdot}-1)c_{k\cdot
}(c_{k\cdot}-1)}{8}Q_{ik}^A
Q^B\\
&&\qquad\quad{} + \sum_{j,l=1}^L \frac{c_{\cdot j}(c_{\cdot j}-1)c_{\cdot
l}(c_{\cdot l}-1)}{8} Q^AQ_{jl}^B\nonumber\\
&&\qquad\quad {}-\frac{\theta_A+\theta_B-c(c-5)+2a+2b-4}{4}\sum_{i=1}^K\sum_{j=1}^L
c_{ij}(c_{ij}-1)Q_i^AQ_j^B\nonumber\\
&&\qquad\quad {}+\sum_{i=1}^K\sum_{j=1}^L \biggl[\frac{c_{i\cdot}(c_{i\cdot
}-1)c_{\cdot j}(c_{\cdot j}-1)}{4}\\
&&\qquad\quad\hspace*{48pt}{} + \frac{c_{ij}(c_{ij}+1-2c_{i\cdot
}+2c_{i\cdot}c_{\cdot j}-2c_{\cdot j})}{2} \nonumber\\
&&\qquad\quad\hspace*{48pt}{}+ c_{ij}b_j(c_{i\cdot
}-1) + c_{ij}a_i(c_{\cdot j}-1) + 2a_i b_j c_{ij} \nonumber\\
&&\qquad\quad\hspace*{48pt}{}+ \frac{\theta
_B}{2}\delta
_{b_j,0}\delta_{c_{\cdot j},1}\delta_{c_{ij},1}(c_{i\cdot} - 1)\\
&&\qquad\quad\hspace*{117pt}{} +
\frac
{\theta_A}{2}\delta_{a_i,0}\delta_{c_{i\cdot},1}\delta
_{c_{ij},1}(c_{\cdot j} -
1) \biggr]Q_i^AQ_j^B\nonumber\\
&&\qquad\quad {}+ \frac{1}{2}\sum_{i=1}^K \biggl[a_i+c_{i\cdot}-1-\frac{\theta_A
}{2}\delta_{c_{i\cdot},2} \biggr]\sum_{j=1}^L c_{ij}(c_{ij}-1)Q_{ii}^A
Q_j^B\nonumber\\
&&\qquad\quad {}+ \frac{1}{2}\sum_{j=1}^L \biggl[b_j+c_{\cdot j}-1-\frac{\theta_B
}{2}\delta_{c_{\cdot j},2} \biggr]\sum_{i=1}^K c_{ij}(c_{ij}-1)Q_i^A
Q_{jj}^B\nonumber\\
&&\qquad\quad {}- \frac{1}{4}\sum_{i=1}^K\sum_{j=1}^L\sum_{k=1}^K
c_{ij}(c_{ij}-1)c_{k\cdot}(c_{k\cdot}-1)Q_{ik}^AQ_j^B\nonumber\\
&&\qquad\quad {}- \frac{1}{4}\sum_{i=1}^K\sum_{j=1}^L\sum_{l=1}^L
c_{ij}(c_{ij}-1)c_{\cdot l}(c_{\cdot l}-1)Q_i^AQ_{jl}^B\nonumber\\
&&\qquad\quad {}+ \frac{1}{8}\sum_{i=1}^K\sum_{j=1}^L\sum_{k=1}^K\sum
_{l=1}^L c_{ij}(c_{ij}-1)c_{kl}(c_{kl}-1)Q_{ik}^AQ_{jl}^B\nonumber
\\
&&\qquad\quad {}- \frac{1}{12}\sum_{i=1}^K\sum_{j=1}^L
c_{ij}(c_{ij}-1)(2c_{ij}-1)Q_{ii}^AQ_{jj}^B.
\end{eqnarray*}

To check the correctness of the above expression, we also solved the recursion
(\ref{secondrec}) numerically for all sample configurations of sizes
$n=10,20$, and $30$ (with $K,L \leq2$), and confirmed that the above
analytic expression agreed in all cases. We also
implemented a \textit{Mathematica} program to solve $q(a,b,c)$ exactly
in the special case $K=L=1$. The program can return series expansions
in terms of $\rho^{-1}$ as $\rho\to\infty$ which are symbolic in
$\theta_A$
and $\theta_B$. We could then compare the first three terms against $q_0$,
$q_1$, and $q_2$, for various sample configurations $(a,b,c)$.
\end{appendix}


\section*{Acknowledgments}
We thank Robert C. Griffiths, Charles H. Langley
and Joshua Paul for useful discussions.

\printaddresses


\begin{thebibliography}{29}

\bibitem[\protect\citeauthoryear{Arratia, Barbour and
  Tavar\'{e}}{2003}]{ABTbook}
\begin{bbook}[author]
\bauthor{\bsnm{Arratia},~\bfnm{A.}\binits{A.}},
  \bauthor{\bsnm{Barbour},~\bfnm{A.~D.}\binits{A.~D.}} \AND
  \bauthor{\bsnm{Tavar\'{e}},~\bfnm{S.}\binits{S.}}
(\byear{2003}).
\btitle{Logarithmic Combinatorial Structures: A Probabilistic Approach}.
\bpublisher{European Mathematical Society Publishing House}, \baddress{Switzerland}.
\bmrnumber{MR2032426}
\end{bbook}
\endbibitem

\bibitem[\protect\citeauthoryear{De Iorio and
  Griffiths}{2004a}]{deigri2004a}
\begin{barticle}[author]
\bauthor{\bsnm{{De Iorio}},~\bfnm{M.}\binits{M.}} \AND
  \bauthor{\bsnm{Griffiths},~\bfnm{R.~C.}\binits{R.~C.}}
(\byear{2004}a).
\btitle{Importance sampling on coalescent histories. {I}}.
\bjournal{Adv. in Appl. Probab.}
\bvolume{36}
\bpages{417--433}.
\bmrnumber{MR2058143}
\end{barticle}
\endbibitem

\bibitem[\protect\citeauthoryear{De Iorio and
  Griffiths}{2004b}]{deigri2004b}
\begin{barticle}[author]
\bauthor{\bsnm{{De Iorio}},~\bfnm{M.}\binits{M.}} \AND
  \bauthor{\bsnm{Griffiths},~\bfnm{R.~C.}\binits{R.~C.}}
(\byear{2004}b).
\btitle{Importance sampling on coalescent histories. {II}}.
\bjournal{Adv. in Appl. Probab.}
\bvolume{36}
\bpages{434--454}.
\bmrnumber{MR2058144}
\end{barticle}
\endbibitem

\bibitem[\protect\citeauthoryear{Ethier and Griffiths}{1990}]{ethgri1990}
\begin{barticle}[author]
\bauthor{\bsnm{Ethier},~\bfnm{S.~N.}\binits{S.~N.}} \AND
  \bauthor{\bsnm{Griffiths},~\bfnm{R.~C.}\binits{R.~C.}}
(\byear{1990}).
\btitle{On the two-locus sampling distribution}.
\bjournal{J. Math. Biol.}
\bvolume{29}
\bpages{131--159}.
\bmrnumber{MR1116000}
\end{barticle}
\endbibitem

\bibitem[\protect\citeauthoryear{Ewens}{1972}]{ewe1972}
\begin{barticle}[author]
\bauthor{\bsnm{Ewens},~\bfnm{W.~J.}\binits{W.~J.}}
(\byear{1972}).
\btitle{The sampling theory of selectively neutral alleles}.
\bjournal{Theor. Popul. Biol.}
\bvolume{3}
\bpages{87--112}.
\bmrnumber{MR0325177}
\end{barticle}
\endbibitem

\bibitem[\protect\citeauthoryear{Fearnhead and Donnelly}{2001}]{feadon2001}
\begin{barticle}[author]
\bauthor{\bsnm{Fearnhead},~\bfnm{P.}\binits{P.}} \AND
  \bauthor{\bsnm{Donnelly},~\bfnm{P.}\binits{P.}}
(\byear{2001}).
\btitle{Estimating recombination rates from population genetic data}.
\bjournal{Genetics}
\bvolume{159}
\bpages{1299--1318}.
\end{barticle}
\endbibitem

\bibitem[\protect\citeauthoryear{Golding}{1984}]{gol1984}
\begin{barticle}[author]
\bauthor{\bsnm{Golding},~\bfnm{G.~B.}\binits{G.~B.}}
(\byear{1984}).
\btitle{The sampling distribution of linkage disequilibrium}.
\bjournal{Genetics}
\bvolume{108}
\bpages{257--274}.
\end{barticle}
\endbibitem

\bibitem[\protect\citeauthoryear{Griffiths}{1981}]{gri1981}
\begin{barticle}[author]
\bauthor{\bsnm{Griffiths},~\bfnm{R.~C.}\binits{R.~C.}}
(\byear{1981}).
\btitle{Neutral two-locus multiple allele models with recombination}.
\bjournal{Theor. Popul. Biol.}
\bvolume{19}
\bpages{169--186}.
\bmrnumber{MR630871}
\end{barticle}
\endbibitem

\bibitem[\protect\citeauthoryear{Griffiths}{1991}]{gri1991}
\begin{bincollection}[author]
\bauthor{\bsnm{Griffiths},~\bfnm{R.~C.}\binits{R.~C.}}
(\byear{1991}).
\btitle{The two-locus ancestral graph}.
In \bbooktitle{Selected Proceedings of the Sheffield Symposium on Applied
  Probability}.
\bseries{IMS Lecture Notes---Monograph Series}
(\beditor{I.~V. Basawa and R.~L. Taylor}, eds.)
\bvolume{18}
\bpages{100--117}.
\bpublisher{IMS}, \baddress{Hayward, CA}.
\bmrnumber{MR1193063}
\end{bincollection}
\endbibitem

\bibitem[\protect\citeauthoryear{Griffiths, Jenkins and
  Song}{2008}]{grietal2008}
\begin{barticle}[author]
\bauthor{\bsnm{Griffiths},~\bfnm{R.~C.}\binits{R.~C.}},
  \bauthor{\bsnm{Jenkins},~\bfnm{P.~A.}\binits{P.~A.}} \AND
  \bauthor{\bsnm{Song},~\bfnm{Y.~S.}\binits{Y.~S.}}
(\byear{2008}).
\btitle{Importance sampling and the two-locus model with subdivided population
  structure}.
\bjournal{Adv. in Appl. Probab.}
\bvolume{40}
\bpages{473--500}.
\bmrnumber{MR2433706}
\end{barticle}
\endbibitem

\bibitem[\protect\citeauthoryear{Griffiths and
  Lessard}{2005}]{Griffiths2005p1491}
\begin{barticle}[author]
\bauthor{\bsnm{Griffiths},~\bfnm{R.~C.}\binits{R.~C.}} \AND
  \bauthor{\bsnm{Lessard},~\bfnm{S.}\binits{S.}}
(\byear{2005}).
\btitle{Ewens' sampling formula and related formulae: Combinatorial proofs,
  extensions to variable population size and applications to ages of alleles}.
\bjournal{Theor. Popul. Biol.}
\bvolume{68}
\bpages{167--177}.
\end{barticle}
\endbibitem

\bibitem[\protect\citeauthoryear{Griffiths and Marjoram}{1996}]{grimar1996}
\begin{barticle}[author]
\bauthor{\bsnm{Griffiths},~\bfnm{R.~C.}\binits{R.~C.}} \AND
  \bauthor{\bsnm{Marjoram},~\bfnm{P.}\binits{P.}}
(\byear{1996}).
\btitle{Ancestral inference from samples of {DNA} sequences with
  recombination}.
\bjournal{J. Comput. Biol.}
\bvolume{3}
\bpages{479--502}.
\end{barticle}
\endbibitem

\bibitem[\protect\citeauthoryear{Hoppe}{1984}]{hop1984}
\begin{barticle}[author]
\bauthor{\bsnm{Hoppe},~\bfnm{F.}\binits{F.}}
(\byear{1984}).
\btitle{P\'{o}lya-like urns and the {E}wens' sampling formula}.
\bjournal{J. Math. Biol.}
\bvolume{20}
\bpages{91--94}.
\bmrnumber{MR758915}
\end{barticle}
\endbibitem

\bibitem[\protect\citeauthoryear{Hudson}{1985}]{hud1985}
\begin{barticle}[author]
\bauthor{\bsnm{Hudson},~\bfnm{R.~R.}\binits{R.~R.}}
(\byear{1985}).
\btitle{The sampling distribution of linkage disequilibrium under an infinite
  allele model without selection}.
\bjournal{Genetics}
\bvolume{109}
\bpages{611--631}.
\end{barticle}
\endbibitem

\bibitem[\protect\citeauthoryear{Hudson}{2001}]{hud2001}
\begin{barticle}[author]
\bauthor{\bsnm{Hudson},~\bfnm{R.~R.}\binits{R.~R.}}
(\byear{2001}).
\btitle{Two-locus sampling distributions and their application}.
\bjournal{Genetics}
\bvolume{159}
\bpages{1805--1817}.
\end{barticle}
\endbibitem

\bibitem[\protect\citeauthoryear{Jenkins and Song}{2009}]{jenson2009}
\begin{barticle}[vtex]
\bauthor{\bsnm{Jenkins},~\bfnm{P.~A.}\binits{P.~A.}} \AND
  \bauthor{\bsnm{Song},~\bfnm{Y.~S.}\binits{Y.~S.}}
(\byear{2009}).
\btitle{Closed-form two-locus sampling distributions: Accuracy and
  universality}.
\bjournal{Genetics}
\bvolume{183}
\bpages{1087--1103}.
\end{barticle}
\endbibitem

\bibitem[\protect\citeauthoryear{Kingman}{1982a}]{kin1982b}
\begin{barticle}[author]
\bauthor{\bsnm{Kingman},~\bfnm{J.~F.~C.}\binits{J.~F.~C.}}
(\byear{1982}a).
\btitle{The coalescent}.
\bjournal{Stochastic Process. Appl.}
\bvolume{13}
\bpages{235--248}.
\bmrnumber{MR671034}
\end{barticle}
\endbibitem

\bibitem[\protect\citeauthoryear{Kingman}{1982b}]{kin1982a}
\begin{barticle}[author]
\bauthor{\bsnm{Kingman},~\bfnm{J.~F.~C.}\binits{J.~F.~C.}}
(\byear{1982}b).
\btitle{On the genealogy of large populations}.
\bjournal{J. Appl. Probab.}
\bvolume{19}
\bpages{27--43}.
\bmrnumber{MR633178}
\end{barticle}
\endbibitem

\bibitem[\protect\citeauthoryear{Kuhner, Yamato and
  Felsenstein}{2000}]{kuhetal2000}
\begin{barticle}[author]
\bauthor{\bsnm{Kuhner},~\bfnm{M.~K.}\binits{M.~K.}},
  \bauthor{\bsnm{Yamato},~\bfnm{J.}\binits{J.}} \AND
  \bauthor{\bsnm{Felsenstein},~\bfnm{J.}\binits{J.}}
(\byear{2000}).
\btitle{Maximum likelihood estimation of recombination rates from population
  data}.
\bjournal{Genetics}
\bvolume{156}
\bpages{1393--1401}.
\end{barticle}
\endbibitem

\bibitem[\protect\citeauthoryear{McVean et al.}{2004}]{McVeanScience04}
\begin{barticle}[author]
\bauthor{\bsnm{McVean},~\bfnm{G.~A.~T.}\binits{G.~A.~T.}},
  \bauthor{\bsnm{Myers},~\bfnm{S.}\binits{S.}},
  \bauthor{\bsnm{Hunt},~\bfnm{S.}\binits{S.}},
  \bauthor{\bsnm{Deloukas},~\bfnm{P.}\binits{P.}},
  \bauthor{\bsnm{Bentley},~\bfnm{D.~R.}\binits{D.~R.}} \AND
  \bauthor{\bsnm{Donnelly},~\bfnm{P.}\binits{P.}}
(\byear{2004}).
\btitle{The fine-scale structure of recombination rate variation in the human
  genome}.
\bjournal{Science}
\bvolume{304}
\bpages{581--584}.
\end{barticle}
\endbibitem

\bibitem[\protect\citeauthoryear{Myers et al.}{2005}]{MyersScience05}
\begin{barticle}[author]
\bauthor{\bsnm{Myers},~\bfnm{S.}\binits{S.}},
  \bauthor{\bsnm{Bottolo},~\bfnm{L.}\binits{L.}},
  \bauthor{\bsnm{Freeman},~\bfnm{C.}\binits{C.}},
  \bauthor{\bsnm{McVean},~\bfnm{G.}\binits{G.}} \AND
  \bauthor{\bsnm{Donnelly},~\bfnm{P.}\binits{P.}}
(\byear{2005}).
\btitle{A fine-scale map of recombination rates and hotspots across the human
  genome}.
\bjournal{Science}
\bvolume{310}
\bpages{321--324}.
\end{barticle}
\endbibitem

\bibitem[\protect\citeauthoryear{Nielsen}{2000}]{nie2000}
\begin{barticle}[author]
\bauthor{\bsnm{Nielsen},~\bfnm{R.}\binits{R.}}
(\byear{2000}).
\btitle{Estimation of population parameters and recombination rates from single
  nucleotide polymorphisms}.
\bjournal{Genetics}
\bvolume{154}
\bpages{931--942}.
\end{barticle}
\endbibitem

\bibitem[\protect\citeauthoryear{Pitman}{1992}]{pitman1992}
\begin{bmisc}[author]
\bauthor{\bsnm{Pitman},~\bfnm{J.}\binits{J.}}
(\byear{1992}).
\btitle{The two-parameter generalization of {E}wens' random partition
  structure}.
Technical Report 345, Dept. Statistics, Univ. California, Berkeley.
\end{bmisc}
\endbibitem

\bibitem[\protect\citeauthoryear{Pitman}{1995}]{pitman1995}
\begin{barticle}[author]
\bauthor{\bsnm{Pitman},~\bfnm{Jim}\binits{J.}}
(\byear{1995}).
\btitle{Exchangeable and partially exchangeable random partitions}.
\bjournal{Probab. Theory Related Fields}
\bvolume{102}
\bpages{145--158}.
\bmrnumber{MR1337249}
\end{barticle}
\endbibitem

\bibitem[\protect\citeauthoryear{Slatkin}{1994}]{SlatkinGenetRes94}
\begin{barticle}[author]
\bauthor{\bsnm{Slatkin},~\bfnm{M.}\binits{M.}}
(\byear{1994}).
\btitle{An exact test for neutrality based on the {E}wens sampling
  distribution}.
\bjournal{Genet. Res.}
\bvolume{64}
\bpages{71--74}.
\end{barticle}
\endbibitem

\bibitem[\protect\citeauthoryear{Slatkin}{1996}]{SlatkinGenetRes96}
\begin{barticle}[vtex]
\bauthor{\bsnm{Slatkin},~\bfnm{M.}\binits{M.}}
(\byear{1996}).
\btitle{A correction to an exact test based on the
Ewens sampling distribution}.
\bjournal{Genet. Res.}
\bvolume{68}
\bpages{259--260}.
\end{barticle}
\endbibitem

\bibitem[\protect\citeauthoryear{Stephens}{2001}]{ste2001}
\begin{bincollection}[author]
\bauthor{\bsnm{Stephens},~\bfnm{M.}\binits{M.}}
(\byear{2001}).
\btitle{Inference under the coalescent}.
In \bbooktitle{Handbook of Statistical Genetics}
(\beditor{D. Balding, M. Bishop and
C. Cannings}, eds.)
\bpages{213--238}.
\bpublisher{Wiley}, \baddress{Chichester, UK}.
\end{bincollection}
\endbibitem

\bibitem[\protect\citeauthoryear{Stephens and Donnelly}{2000}]{stedon2000}
\begin{barticle}[author]
\bauthor{\bsnm{Stephens},~\bfnm{M.}\binits{M.}} \AND
  \bauthor{\bsnm{Donnelly},~\bfnm{P.}\binits{P.}}
(\byear{2000}).
\btitle{Inference in molecular population genetics}.
\bjournal{J. R. Stat. Soc. Ser. B Stat.
Methodol.}
\bvolume{62}
\bpages{605--655}.
\bmrnumber{MR1796282}
\end{barticle}
\endbibitem

\bibitem[\protect\citeauthoryear{Wang and Rannala}{2008}]{wanran2008}
\begin{barticle}[author]
\bauthor{\bsnm{Wang},~\bfnm{Y.}\binits{Y.}} \AND
  \bauthor{\bsnm{Rannala},~\bfnm{B.}\binits{B.}}
(\byear{2008}).
\btitle{Bayesian inference of fine-scale recombination rates using population
  genomic data}.
\bjournal{Philos. Trans. R. Soc.}
\bvolume{363}
\bpages{3921--3930}.
\end{barticle}
\endbibitem

\bibitem[\protect\citeauthoryear{Watterson}{1977}]{Watterson1977}
\begin{barticle}[author]
\bauthor{\bsnm{Watterson},~\bfnm{G.~A.}\binits{G.~A.}}
(\byear{1977}).
\btitle{Heterosis or neutrality?}
\bjournal{Genetics}
\bvolume{85}
\bpages{789--814}.
\bmrnumber{MR0504021}
\end{barticle}
\endbibitem


\end{thebibliography}
\end{document}